\documentclass{m2an}

\newcommand{\rd}{\mathrm{d}}

\usepackage[colorlinks,linkcolor=blue,anchorcolor=blue,citecolor=red]{hyperref}

\newtheorem{theorem}{Theorem}

\newtheorem{lemma}[theorem]{Lemma}
\newtheorem{remark}{Remark}

\usepackage{bm}
\usepackage{float}
\usepackage{multirow}
\usepackage{upgreek}
\usepackage{amsmath}
\usepackage{amsfonts}
\usepackage{amssymb}
\usepackage{graphicx}
\usepackage{mathrsfs}
\usepackage{tabularx}
\usepackage{array}
\usepackage{cite}
\usepackage{amsthm}

\begin{document}
\title{Error estimates of asymptotic-preserving neural networks in approximating stochastic linearized Boltzmann equation}
\author{Jiayu Wan}\address{The Chinese University of Hong Kong, Hong Kong}
\author{Liu Liu}\address{The Chinese University of Hong Kong, Hong Kong}

\begin{abstract}
In this paper, we construct an asymptotic-preserving neural networks (APNNs) \cite{JMW2023} for the linearized Boltzmann equation in the acoustic scaling and with uncertain parameters. Utilizing the micro-macro decomposition, we design the loss function based on the stochastic-Galerkin system conducted from the micro-macro equations. Rigorous analysis is provided to show the capability of neural networks in approximating solutions near the global Maxwellian. By employing hypocoercivity techniques, we demonstrate two key results: (i) the existence of APNNs when the loss function approaches zero, and (ii) the convergence of the APNN's approximated solution as the loss tends to zero, with the error exhibiting an exponential decay in time.
\end{abstract}

\subjclass{35Q20, 68T07, 82C40, 65F99}
\keywords{linearized Boltzmann equation \and uncertainty quantification \and deep learning \and asymptotic-preserving neural networks \and hypocoercivity}
\date{}

\maketitle
Kinetic equations have been widely used in many areas such as rarefied gas, plasma physics, astrophysics, semiconductor device modeling, and social and biological sciences \cite{Semi-Book, Villani02}. They describe the non-equilibrium dynamics of a system composed of a large number of particles and bridge atomistic and continuum models in the hierarchy of multiscale modeling. The Boltzmann-type equation, as one of the most representative models in kinetic theory, provides a power tool to describe molecular gas dynamics, radiative transfer, plasma physics, and polymer flow \cite{A15}. They have significant impacts in designing, optimization, control, and inverse problems. For example, it can be used in the design of semiconductor devices, topology optimization of gas flow channels, or risk management in quantitative finance \cite{CS21}. Many of these applications often require finding unknown or optimal parameters in the Boltzmann-type equations or mean-field models \cite{AHP15, Caflisch, CFP13, Cheng11}. 

In addition, kinetic equations typically involve various sources of uncertainty, such as modeling errors, imprecise measurements, and uncertain initial conditions. In particular, the collision kernel in the Boltzmann equation governs the transition rates during particle collisions. Calculating this collision kernel from first principles is highly complex, and in practice, heuristic approximations or empirical data are often used, inevitably introducing uncertainties. Additionally, uncertainties may stem from inaccurate measurements of initial or boundary conditions, as well as from source terms, further compounding the uncertainties in the model. As a result, addressing uncertainty quantification (UQ) becomes essential for evaluating, validating, and improving the underlying models, underscoring our project's significance. For numerical studies of the Boltzmann equation and other kinetic models with or without randomness, we refer readers to works such as \cite{HPY, JXZ15, Jin-ICM, PZ2020} and \cite{pareschi, MP06, MPR13, HQ20, HQY21}.
Among the various numerical approaches, the generalized polynomial chaos (gPC)-based stochastic Galerkin (SG) method and its variations have been widely adopted, demonstrating success in a range of applications \cite{Xiu}. Beyond numerical simulations, theoretical studies have established the stability and convergence of these methods. Spectral convergence for the gPC-based SG method was demonstrated in \cite{HJ16, LQ2024-1, LQ2024-2}, while \cite{LiuJin2018} and \cite{LiuDaus2019} introduced a robust framework based on hypocoercivity to perform local sensitivity analysis for a class of multiscale, inhomogeneous kinetic equations with random uncertainties—approximated using the gPC-based SG method. For further reference, we point readers to the recent collection \cite{JinPareschi} and the survey \cite{parUQ}.

Modeling and predicting the evolution of multiscale systems such as the Boltzmann-type equations have always been challenging, which often requires sophisticated knowledge of numerical methods and labor-intensive implementation, in addition to the prohibitive costs due to the well-known curse of dimensionality. The issue of high dimensionality becomes even more overwhelming when uncertainties are considered, making traditional numerical approaches unfavorable. This motivates researchers to develop data-driven models and methods \cite{WE} in the recent decade.

Machine Learning, or Deep neural networks (DNNs) in particular, have gained increasing interest in approximating the solutions of partial differential equations (PDE) due to their universal approximation property and ability to handle high-dimensional problems. Many studies have been done to employ DNN for solving deterministic or parameterized PDEs, and have shown remarkable promise in various applications. Among those, \textit{Physics Informed Neural Network} (PINN) \cite{MPG2019} is one of the most famous approaches which incorporates the available physical laws and limited data, such as boundary or initial conditions and the source term, to approximate the solution of the underlying PDE. The idea has been successfully applied in the simulation of many forward and inverse problems \cite{CLKD2020,LMK2021,MJK2020}.

However, if one directly applies the standard PINNs when dealing with multiscale kinetic models, it may lead to incorrect inferences and predictions. This is due to the presence of small scales needing to be enforced consistently during the learning process, but a standard PINN formulation only captures the solution at the leading order of the Knudsen number, thereby losing accuracy in the asymptotic limit regimes. To overcome this difficulty, the authors in \cite{JMW2023} propose the \textit{Asymptotic-Preserving Neural Networks} (APNNs) to enhance the performance of standard PINN to solve multiscale linear transport equations. The idea of APNN is inspired by the traditional AP schemes(\cite{Jin2022}), which preserve the asymptotic transition from one scale to another at the discrete level and capture the limiting macroscopic behavior of the solution when the scaling parameter approaches zero. Recently, an APNN-based method is performed to solve the linear semiconductor Boltzmann equation with good accuracy \cite{LiuZhu2024}. APNN is also applied to study hyperbolic-type linear kinetic equations with multiple scales, for example, see \cite{BLPZ2022}. 

Although some numerical experiments have been conducted to validate the efficiency and accuracy of APNN-based methods, rigorous analysis of the convergence of these methods is still limited. For literature in this direction, we refer to \cite{Yangxu2024}, in which the authors study Boltzmann equations with linear collision kernels and present a formal proof of the convergence of the APNN solutions to the real solutions in the standard $L^2$ space, as the defined APNN loss approaches zero. In this paper, we carry out rigorous analysis of APNN for linearized Boltzmann equation perturbed around some global Maxwellian. We will perform error estimates for the method and derive some convergence results. Our innovations are twofold. First, we introduce uncertainties into our system and adopt the SG method to show the convergence of the approximated SG coefficients in some weighted normed space, as inspired by \cite{LiuJin2018}. Second, we conduct hypocoercivity analysis \cite{briant2015, MN2006} to the APNN system to derive convergence in $H^1$ space (not $L^2$ space) with errors decaying exponentially in time, thus providing a more accurate rate of convergence compared to \cite{Yangxu2024}. Our analysis can be outlined as follows. First, we perform micro-macro decomposition for the linearized Boltzmann equation and define our APNN loss function based on the micro-marco system. We then have two main theorems to prove: the first theorem suggests that there exists neural networks which lead to arbitrarily small APNN loss, and the proof is based on the Universal Approximation Theorem(UAT) of neural networks. The second theorem states that the errors of the APNN-approximated solutions tend to zero as the APNN loss approaches zero, with an exponential decay in time. The proof is inspired by \cite{briant2015}, which constructs a Lyapunov functional that is equivalent to the standard Sobolev norm $H^1$. We then estimate the time evolution for each term involved in the functional and keep track of the terms that serve as components of the APNN loss function. 

The rest of the paper is organized as follows. In Section 1, we review some important concepts in kinetic theory which serve as the cornerstone for our analysis. This includes a short introduction to linearized Boltzmann equation with uncertainty, a complete list of hypocoercivity assumptions, and a concise review of the stochastic Galerkin method. In Section 2, we derive a coupled system in the micro-macro decomposition framework with a formal analysis of the asymptotic limit. We then define our APNN loss function based on the decomposition and generalize the result to the stochastic Galerkin setting. In Section 3, we prove our main results, which are the two theorems mentioned in the previous paragraph.

\section{Preliminaries}
\label{preliminary}

\subsection{\textbf{Introduction to Boltzmann equation with uncertainties}} Consider the initial value problem for the Boltzmann equation

\begin{equation}\label{Boltzmann equation}
\begin{cases}
    \ \partial_{t}f + \frac{1}{\varepsilon^{\alpha}}v\cdot \nabla_{x}f = \frac{1}{\varepsilon^{1+\alpha}}Q(f,f) \\[4pt]
    \ f(0,x,v,z)= f_{I}(x,v,z), ~~~~x\in \mathbb{T}^{d}, v \in \mathbb{R}^{d}, z \in I_{z} \subset \mathbb{R}
\end{cases}
\end{equation}

where $f=f(t,x,v,z)$ represents the particle density distribution in the phase space that depends on time $t$, particle position $x$, particle velocity $v$ and a random variable $z$. The number $d \geq 1$ denotes the dimension of the spatial and velocity spaces, and $z$ is a random variable that lies in domain $I_{z} \subset \mathbb{R}$ with compact support, which is used to account for the random uncertainties of inputs. The operator $Q$ is quadratic and models the binary collisional interactions between particles. The parameter $\varepsilon$ is the dimensionless Knudsen number, the ratio of particle mean free path over the domain size. The choice $\alpha=1$ is referred to the incompressible Navier-Stokes scaling, while $\alpha=0$ corresponds to the Euler(acoustic) scaling. Periodic boundary condition for the spatial domain $\Omega=\mathbb{T}^{d}$ is assumed. 

We consider random collision kernels, hence the operator $Q$ is defined by

\begin{equation}\label{collision-kernel}
    Q(f,g)=\int_{\mathbb{R}^d \times \mathbb{S}^{d-1}} B(|v-v_{*}|,\cos{\theta},z)(f^{'}g^{'}_{*}+f^{'}_{*}g^{'}-fg_{*}-f_{*}g) ~dv_{*}d\sigma
\end{equation}

\noindent where we adopt the notations $f^{'}=f(v^{'}),f_{*}=f(v_{*}), f^{'}_{*}=f(v^{'}_{*})$ and similar for $g$. $v^{'}$ and $v_{*}^{'}$ are the post-collisional velocities of particles depending on the pre-collisional velocities $v$ and $v_{*}$, which are given by:

\begin{equation*}
    v^{'}=\frac{v+v_{*}}{2}+\frac{|v-v_{*}|}{2}\sigma, ~~~~~~v^{'}_{*}=\frac{v+v_{*}}{2}-\frac{|v-v_{*}|}{2}\sigma
\end{equation*}

\noindent $\theta \in [0,\pi]$ is the deviation angle between $v^{'}-v^{'}_{*}$ and $v-v_{*}$. The collision kernel $B=B(|v-v^{*}|,\cos{\theta},z)$ is a non-negative function determined by physics, and it is assumed to depend on the random variable $z \in I_{z}$ in our setting. The Boltzmann collision operator conserves mass, momentum and energy, namely
\begin{equation}\label{conservation_laws}
    \int_{\mathbb{R}^d} Q(f,f)\phi(v) ~dv=0, ~~~\phi(v)=1,v,|v|^2
\end{equation}

\noindent Moreover, we have the celebrated Boltzmann's H-theorem, 

\begin{equation}\label{H_theorem}
    \int_{\mathbb{R}^d} Q(f,f)\log{(f)} ~dv \leq 0, 
\end{equation}
such that 
\begin{equation*}
    \int_{\mathbb{R}^d} Q(f,f)\log{(f)} ~dv =0 ~\Leftrightarrow~ Q(f,f)=0 ~\Leftrightarrow~ f=\widetilde{M},  
\end{equation*}
where $\widetilde{M}$ is the \textit{local} equilibrium state: 
\begin{equation}\label{local_Maxwellian}
    \widetilde{M}=\frac{\rho}{(2\pi T)^{\frac{d}{2}}}\exp{(-\frac{|v-u|^2}{2T})}, 
\end{equation}
 with 
\begin{equation*}
    \begin{aligned}
        \rho(x,t)=& \int_{\mathbb{T}^d \times \mathbb{R}^d} f ~dv, \\
        u(x,t)=& \frac{1}{\rho}\int_{\mathbb{T}^d \times \mathbb{R}^d} fv ~dv, \\
        T(x,t)=& \frac{1}{d\rho}\int_{\mathbb{T}^d \times \mathbb{R}^d} f|v-u|^2 ~dv, 
    \end{aligned}
\end{equation*}
 corresponding to the density, mean velocity and temperature of the gas respectively, which are all determined by the initial datum due to the conservation properties. The \textit{global} equilibrium is the unique stationary solution to \eqref{Boltzmann equation} and is given by
\begin{equation}\label{global_Maxwellian}
    \mathcal{M}(v) = \frac{1}{(2\pi)^{\frac{d}{2}}}e^{-\frac{|v|^2}{2}}. 
\end{equation}

We consider a linearization around the \textit{global} equilibrium for the solution: 
\begin{equation}\label{linear_perturbation}
    f=\mathcal{M}+\varepsilon M h, 
\end{equation}
where $M=\sqrt{\mathcal{M}}$. Inserting \eqref{linear_perturbation} into the Boltzmann equation \eqref{Boltzmann equation}, the fluctuation $h$ satisfies
\begin{equation}\label{linearized_Boltzmann_equation}
    \begin{cases}
    \displaystyle \partial_{t}h + \frac{1}{\varepsilon^{\alpha}}v\cdot \nabla_{x}h = \frac{1}{\varepsilon^{1+\alpha}}\mathcal{L}(h)+ \frac{1}{\varepsilon^{\alpha}}\mathcal{F}(h,h) \\[4pt]
    \ h(0,x,v,z)= h_{I},
\end{cases}
\end{equation}
where the linearized collision operator $\mathcal{L}$ is defined as 
\begin{equation}\label{linearized_operator}
    \begin{aligned}
        \mathcal{L}(h)=& M^{-1}\left(Q(Mh,\mathcal{M})+Q(\mathcal{M},Mh)\right) \\
        =& M\int_{\mathbb{R}^d \times \mathbb{S}^{d-1}} B(|v-v_{*}|,\cos{\theta},z)\mathcal{M}_{*}\left(\frac{h^{'}_{*}}{M^{'}_{*}} + \frac{h^{'}}{M^{'}} - \frac{h_{*}}{M_{*}} - \frac{h}{M}\right) dv_{*}d\sigma,
    \end{aligned}
\end{equation}
with the nonlinear operator $\mathcal{F}$ given by 
\begin{equation}\label{nonlinear_operator}
    \begin{aligned}
        \mathcal{F}(h,h)=& M^{-1}[Q(Mh,Mh)+Q(Mh,Mh)] \\
        =& \int_{\mathbb{R}^d \times \mathbb{S}^{d-1}} B(|v-v_{*}|,\cos{\theta},z)M_{*}\left(h^{'}_{*}h^{'}-h_{*}h\right)~dv_{*}d\sigma
    \end{aligned}
\end{equation}

In this paper, we will focus mainly on the acoustic scaling($\alpha=0$) in spatial dimension three($d=3$), and we will ignore the influence of the nonlinear operator $\mathcal{F}$ by setting $\mathcal{F}=0$ \cite{MN2006}. Then \eqref{linearized_Boltzmann_equation} becomes 

\begin{equation}\label{fully_linearized_Boltzmann_with_acoustic_scaling}
    \begin{cases}
    \displaystyle \partial_{t}h + v\cdot \nabla_{x}h = \frac{1}{\varepsilon}\mathcal{L}(h) \\
    \ h(0,x,v,z)= h_{I},
    \end{cases}
\end{equation}

The linear operator $\mathcal{L}$ is self-adjoint on the space $L_{v}^{2}$, and it can be split as
\begin{equation*}
    \mathcal{L}(h) = K(h) - \Lambda(h)
\end{equation*}
such that 
\begin{equation*}
    K(h) = L^{+}(h) - L^{*}(h), \qquad
    L^{*}(h) = M\left[(hM)* \Upphi\right], 
\end{equation*}
with 
\begin{equation*}
    L^{+}(h) = \int_{\mathbb{R}^3 \times \mathbb{S}^{2}} B(|v-v_{*}|,\cos{\theta},z)\left(h^{'}M^{'}_{*}+h_{*}^{'}M^{'}\right)M_{*} ~dv_{*}d\sigma. 
\end{equation*}
In addition, $\Lambda(h) = \nu(v, z)h$, with the collision frequency 
\begin{equation}\label{collision_frequency}
    \nu(v,z)=\int_{\mathbb{R}^3 \times \mathbb{S}^{2}} B(|v-v_{*}|,\cos{\theta},z)\mathcal{M}_{*}~dv_{*}d\sigma = (\Upphi*\mathcal{M})(v). 
\end{equation}

\textbf{Assumptions on the collision kernel:}

We will make the same assumptions on the collision kernel as in \cite{LiuDaus2019}. In particular, we will consider hard potentials with $B$ satisfying Grad's angular cutoff, that is, 
\begin{equation}\label{assumptions on collision kernel}
    \begin{aligned}
        & B(|v-v_{*}|,\cos{\theta},z)=\Upphi(|v-v_{*}|)b(\cos{\theta},z), \qquad \Upphi(|v-v_{*}|)=C|v-v_{*}|^{\gamma},~ \gamma \in [0,1],~ C>0, \\
        & \forall \eta \in [-1,1], \qquad |b(\eta,z)| \leq C_{b},~ |\partial_{\eta}b(\eta,z)| \leq \tilde{C_{b}},~ |\partial_{z}^{k}b(\eta,z)| \leq C_{b}^{*}, ~~~\forall ~ 0 \leq k \leq r.
    \end{aligned}
\end{equation}
We assume that $b$ is linear in $z$, namely
\begin{equation}\label{b linear in z}
    b(\cos{\theta},z)=b_{0}(\cos{\theta}) + b_{1}(\cos{\theta})z, ~~~|z| \leq C_{z},
\end{equation}
Furthermore, we presume that
\begin{equation}\label{assumption on b0}
    b_{0}(\cos{\theta}) \geq (2^q+2)|b_{1}(\cos{\theta})|C_{z} +D(\cos{\theta}), 
\end{equation}
\noindent where $D \in C^{1}([-1,1])$ such that $|D(\eta)|,\, |D^{'}(\eta)| \leq C_{D}$, $\forall \eta \in [-1,1]$. The importance of the linearity assumption \eqref{b linear in z} will become clear when we introduce the stochastic Galerkin method in section 1.3. \eqref{assumption on b0} is a technical assumption to ensure that our hypocoercivity analysis can be generalized to the stochastic setting (see Section 3 for more details). 

\bigskip

\subsection{\textbf{Hypocoercivity assumptions.} }

Our error estimate is based on hypocoercivity analysis, which relies on some assumptions on the collision kernel. In this subsection, we will summarize these assumptions and we comment that the linearized Boltzmann operator $\mathcal{L}$ satisfies all these assumptions \cite{MN2006}. 

\bigskip

\noindent \textbf{H1: $\Lambda$ satisfies coercivity conditions.} $\mathcal{L}: L^{2}_{x,v}=L^2(\mathbb{T}^3 \times \mathbb{R}^3)$ is closed, self-adjoint on $L_{v}^{2}$ and local in $x$. $\mathcal{L}$ has the form $\mathcal{L}=K - \Lambda$. There is a norm $||\cdot||_{\Lambda_{v}}$ on $L^{2}_{v}$ given by

\begin{equation}\label{def_Lambda_norm}
    ||h||_{\Lambda_{v}}= ||h(1+|v|)^{\frac{\gamma}{2}}||_{L^2_{v}}
\end{equation}

\noindent such that $\forall h \in L^{2}_{v}$, $\Lambda$ satisfies the coercivity condition: 
\begin{equation}\label{H_1_1}
    \nu_{0}^{\Lambda} ||h||_{L_{v}^{2}}^{2} \leq \nu_{1}^{\Lambda} ||h||_{\Lambda_{v}}^{2} \leq \langle \Lambda(h),h \rangle _{L_{v}^{2}}\leq \nu_{2}^{\Lambda} ||h||_{\Lambda_{v}}^{2},
\end{equation}

\noindent and $\forall h\in H^{1}_{v}$, 
\begin{equation}\label{H_1_2}
    \langle \nabla_{v}\Lambda(h), \nabla_{v}h \rangle_{L_{v}^{2}} \geq \nu_{3}^{\Lambda} ||\nabla_{v}h||_{\Lambda_{v}}^{2} - \nu_{4}^{\Lambda} ||h||_{\Lambda_{v}}^{2},
\end{equation}

\noindent where $(\nu_{s}^{\Lambda})_{0 \leq s \leq 4} >0$ are constants depending on the operators and the velocity space. We also assume that $\forall h,g \in L^{2}_{v}$, 
\begin{equation}\label{H_1_3}
    \langle \mathcal{L}(h),g \rangle_{L_{v}^{2}} \leq C^{\mathcal{L}}||h||_{\Lambda_{v}}^{2} ||g||_{\Lambda_{v}}^{2}
\end{equation}

\noindent We define a norm $||\cdot||_{\Lambda}$ on $L^{2}_{x,v}$ by 
\begin{equation}\label{def_total_Lambda_norm}
    ||\cdot||_{\Lambda}^2 = \int_{\mathbb{T}^3}||\cdot||_{\Lambda_{v}}^2 ~dx. 
\end{equation}

\bigskip

\noindent \textbf{H2: $K$ has a regularizing effect.} For any $ \delta>0$, there exists some explicit constant $C(\delta)>0$ such that $\forall h \in H^{1}_{v}$, 
\begin{equation}\label{H_2}
    \langle \nabla_{v} K(h),\nabla_{v}h \rangle_{L^2_{v}} \leq C(\delta)||h||_{L_{v}^{2}}^{2} + \delta||\nabla_{v}h||_{L_{v}^{2}}^{2}
\end{equation}

\noindent In fact, it is shown in \cite{MN2006} that $ ||\nabla_{v} K(h)||^2_{L_{v}^2} \leq C(\delta)||h||_{L_{v}^{2}}^{2} + \delta||\nabla_{v}h||_{L_{v}^{2}}^{2}$.

\bigskip

\noindent \textbf{H3: $\mathcal{L}$ has a finite dimensional kernel}. One has
$$ N(\mathcal{L}) = Span\{ \varphi_1,...,\varphi_N \}, $$
such that $\{ \varphi_{i}\}_{1 \leq i \leq N}$ forms an orthonormal set with $\varphi_{i}(v) = P_{i}(v)e^{-\frac{|v|^{2}}{4}}$,  where $P_{i}$ is a polynomial. 

We denote by $\pi_{\mathcal{L}}$ the orthogonal projection in $L_{v}^{2}$ onto $N(\mathcal{L})$. Then $\forall h \in L_{v}^{2}$, 

\begin{equation}\label{ortho_projection}
    \pi_{\mathcal{L}}(h)=\sum_{i=1}^{N} \left(\int_{\mathbb{R}^3} h\varphi_{i} ~dv\right) \varphi_{i}. 
\end{equation}

\noindent Here $\mathcal{L}$ has the local coercivity property: there exists $\lambda >0$ such that $\forall h\in L_{v}^{2}$, 

\begin{equation}\label{H_3}
    \langle \mathcal{L}(h),h \rangle_{L_{v}^{2}} \leq -\lambda ||h^{\perp}||_{\Lambda_{v}}^{2}, 
\end{equation}

\noindent where $h^{\perp}=h - \pi_{\mathcal{L}}(h)$ represents the microscopic part of $h$. We summarize here some key facts about the fluid projection $\pi_{\mathcal{L}}$, which will be used later. 

If $h \in H^{1}_{x,v}$, write $\pi_{\mathcal{L}}(h) = \sum h_i \varphi_i$ as in \eqref{ortho_projection}, then $\partial_{v_i}\varphi_j$ is still of the form $P(v)e^{-\frac{|v|^{2}}{4}}$ for a polynomial $P$, whose norm in $L_{v}^{2}$ is finite. We let $M=\max \{ ||\nabla_{v}\varphi_{i}||^2_{L^2_{v}} \}_{1 \leq i \leq N}$, then $||\nabla_{v}\pi_{\mathcal{L}}(h)||_{L_{x,v}^{2}}^{2} \leq M ||\pi_{\mathcal{L}}(h)||_{L_{x,v}^{2}}^{2}$ by Cauchy-Schwarz and the orthonormality of $\{ \varphi_{i} \}$. Similar results can be proved for $||v \cdot \nabla_{v}\pi_{\mathcal{L}}(h)||_{L_{x,v}^{2}}^{2}$ and $ ||\nabla_{v}\pi_{\mathcal{L}}(v \cdot \bm{h})||_{L_{x,v}^{2}}^{2}$ by choosing $M$ to be $\max \{ ||v \cdot \nabla_{v}\varphi_{i}||^2_{L^2_{v}} \}$ and $\max \{ || \nabla_{v}v_{i}\varphi_{j}||^2_{L^2_{v}} \}$ respectively. Hence there exists $C_{\pi 1}>0$ such that $\forall h \in H^{1}_{x,v}$, 
 
 \begin{equation}\label{H_3_1}
    \begin{aligned}
        \ & ||\nabla_{v}\pi_{\mathcal{L}}(h)||_{L_{x,v}^{2}}^{2} \leq C_{\pi 1} ||\pi_{\mathcal{L}}(h)||_{L_{x,v}^{2}}^{2} \\
        \ & ||v \cdot \nabla_{v}\pi_{\mathcal{L}}(h)||_{L_{x,v}^{2}}^{2} \leq C_{\pi 1} ||\pi_{\mathcal{L}}(h)||_{L_{x,v}^{2}}^{2}\\
        \ & ||\nabla_{v}\pi_{\mathcal{L}}(v \cdot \bm{h})||_{L_{x,v}^{2}}^{2} \leq C_{\pi 1} ||\pi_{\mathcal{L}}(\bm{h})||_{L_{x,v}^{2}}^{2}
    \end{aligned}
 \end{equation}

 \noindent Moreover, setting $M=\max \{ ||(1+|v|)^{\frac{\gamma}{2}}\varphi_{i}||^2_{L^2_{v}} \} $ and $M=\max \{ ||(1+|v|)^{\frac{\gamma}{2}}\nabla_{v}\varphi_{i}||^2_{L^2_{v}} \} $ respectively, we can prove similar results for $||\pi_{\mathcal{L}}(h)||^2_{\Lambda}$ and $||\nabla_{v}\pi_{\mathcal{L}}(h)||^2_{\Lambda}$. Hence there exists $C_{\pi} >0 $ such that $\forall h \in H_{x,v}^{1}$, 

 \begin{equation}\label{H_3_2}
    \begin{aligned}
        \ & ||\pi_{\mathcal{L}}(h)||_{\Lambda}^{2} \leq C_{\pi}||\pi_{\mathcal{L}}(h)||_{L_{x,v}^{2}}^{2} \leq C_{\pi}||h||_{L_{x,v}^{2}}^{2} \\
        \ & ||\nabla_{v}\pi_{\mathcal{L}}(h)||_{\Lambda}^{2} \leq C_{\pi}||\pi_{\mathcal{L}}(h)||_{L_{x,v}^{2}}^{2} \leq C_{\pi}||h||_{L_{x,v}^{2}}^{2}
    \end{aligned}
 \end{equation}

Inequality \eqref{H_3_2}, combining with \textbf{H1}, shows that the $\Lambda$-norm and the standard $L_{x,v}^{2}$-norm are equivalent on the fluid part of $\mathcal{L}$. 
Finally, for all $h \in N(\frac{1}{\varepsilon}\mathcal{L} - v\cdot \nabla_{x})^{\perp}$, the Poincare inequality on the torus gives that for some $C_{p}>0$, 
 \begin{equation}\label{Poincare inequality}
     ||\pi_{\mathcal{L}}(h)||^2_{L_{x,v}^{2}} \leq C_{p}||\nabla_{x}\pi_{\mathcal{L}}(h)||^2_{L_{x,v}^{2}} \leq C_{p}||\nabla_{x}h||^2_{L_{x,v}^{2}}. 
 \end{equation}

\bigskip

 \subsection{\textbf{Stochastic Galerkin (SG) method}}
 \label{SG method}
In this subsection, we discuss how to manage the random variable $z$ existing in the perturbed Boltzmann equation \eqref{fully_linearized_Boltzmann_with_acoustic_scaling}. The main idea involves applying the stochastic Galerkin (SG) method to eliminate the randomness, thereby transforming the equation into a system of equations containing only deterministic coefficients. Throughout the paper, we assume that $z \in I_z \subset \mathbb{R}$ with compact support. 

We define the space 
$$\mathbb{P}^K := \text{Span} \left\{ \phi_{i}(z) \, \Big| \, 1 \leq i \leq K \right\}$$
equipped with the inner product with respect to the probability density function $\pi(z)$ in $z$(which is given \textit{a priori}):
\begin{equation*}
	\langle f(t,x,v,\cdot), \, g(t,x,v,\cdot) \rangle_{I_{z}} = \int_{I_{z}}  f(t,x,v,z)g(t,x,v,z)\pi(z) ~dz,
\end{equation*}
where $\{ \phi_{i}(z)\}_{i = 1}^{K}$ is an orthonormal gPC basis function, i.e.,
\begin{equation}\label{orthonormal-relation}
    \int_{I_z} \phi_{i}(z) \phi_{j}(z) \pi(z) \,\rd z = \delta_{ij}, \quad 1 \leq i, j \leq K.
\end{equation}

For \eqref{fully_linearized_Boltzmann_with_acoustic_scaling}, we employ the SG method and expand $h(t,x,v,z)$ as
\begin{equation}\label{h-expansion}
    h(t,x,v,z) \approx h_{K}(t,x,v,z) := \sum_{i=1}^{K} h^{i}(t,x,v) \phi_{i}(z), 
\end{equation}
with
\begin{equation}\label{hk}
    h^{i}(t,x,v)= \int_{I_z} h(t,x,v,z)\phi_{i}(z)\pi(z)dz.
\end{equation}

\noindent We perform the Galerkin projection of \eqref{fully_linearized_Boltzmann_with_acoustic_scaling} onto $\mathbb{P}^{K}$ to obtain a system of equations: 
\begin{equation}\label{SG_linearized_Boltzmann_equation}
    \begin{cases}
    \ \partial_{t}h^i + v\cdot \nabla_{x}h^i = \frac{1}{\varepsilon}\mathcal{L}_{i}(h_K), \\[2pt]
    \ h_i(0,x,v)= h_{i}^{I}(x,v), 
\end{cases}
\end{equation}

\noindent for $1 \leq i \leq K$, with periodic boundary conditions and initial data for $h_i$ given by 
\begin{equation*}
    h_{i}^{I}(x,v) := \int_{I_z} h_{I}(x,v,z)\phi_{i}(z)\pi(z) ~dz. 
\end{equation*}

\noindent The linear term $\mathcal{L}_{i}(h_K)$ is defined by

\begin{equation}\label{SG_linear_term}
\begin{aligned}
    \mathcal{L}_{i}(h_K)=& \langle \mathcal{L}(h_K), \phi_{i} \rangle_{I_z} 
    = \sum_{k=1}^{K} \langle \mathcal{L}(h_{k}\phi_{k}), \phi_{i} \rangle_{I_z} \\
    =& \sum_{k=1}^{K} M\int_{\mathbb{R}^3 \times \mathbb{S}^2} S_{ik}\mathcal{M}_{*}\left(\frac{h_{*}^{k,'}}{M^{'}_{*}} + \frac{h^{k,'}}{M^{'}} - \frac{h_{*}^{k}}{M_{*}} - \frac{h^{k}}{M}\right)dv_{*}d\sigma \\
    \ = &\sum_{k=1}^{K}\mathcal{L}_{ik}(h^k)
\end{aligned}
\end{equation}
with
\begin{equation}\label{S_ik}
    S_{ik}(|v-v_{*}|,\cos{\theta}) = \int_{I_z} B(|v-v_{*}|,\cos{\theta},z)\phi_{i}(z)\phi_{k}(z)\pi(z) ~dz.
\end{equation}

The linearity assumption \eqref{b linear in z} guarantees that for each fixed $i$, $S_{ik}$ is nonzero for only three choices of $k: i-1, i, i+1$, which will simplify our error analysis in Section 3. We hence define 

\begin{equation}\label{def of chi}
    \chi_{ik} = 
    \begin{cases}
        \ 0, ~~S_{ik}=0 \\
        \ 1, ~~S_{ik} \neq 0
    \end{cases}
\end{equation}

Similar to \cite{LiuJin2018}, we make the technical assumption 
\begin{equation*}
    ||\phi_{i}||_{L^{\infty}} \leq Ci^{p}, \qquad \forall 1 \leq i \leq K, 
\end{equation*}
\noindent with parameter $p \geq 0$. Let $q>p+2$, we then define the energy $E^{K}_{t}$ by 

\begin{equation}\label{energy norm def}
    E^{K}_{t}(\bm{h}) = \sum_{i=1}^{K} ||i^q h^{i}(t,\cdot,\cdot)||_{H_{x,v}^{1}}^{2}, 
\end{equation}
\noindent where $\bm{h} = (h^1,h^2,...,h^K)^{T}$.

\section{APNN framework}
\label{APNN frame}

\subsection{\textbf{Micro-macro decomposition}}

In this subsection, we derive the micro-macro decomposition for the linearized Boltzmann equation perturbed around the global equilibrium $\mathcal{M}(v)$. As before, we assume acoustic scaling and suppose first that the model is deterministic with spatial dimension $d=3$. We recall here that the equation is given by 

\begin{equation}\label{linearized_Boltzmann_equation_Euler_scaling}
    \begin{cases}
    \ \partial_{t}h + v\cdot \nabla_{x}h = \frac{1}{\varepsilon}\mathcal{L}(h),  \\
    \ h(0,x,v)= h_{I}. 
\end{cases}
\end{equation}

Denote $M=\sqrt{\mathcal{M}}$. By \textbf{H3}, the null space of $\mathcal{L}$ is finite dimensional and it is well known that $N(\mathcal{L})$ has an orthonormal basis given by
\begin{equation*}
    N(\mathcal{L})=\text{Span} \{ \varphi_{0}M,\varphi_{1}M, \varphi_{2}M, \varphi_{3}M, \varphi_{4}M \}
\end{equation*}
 where 
\begin{equation}\label{o.n. basis expression for L}
    \begin{cases}
    \ \varphi_{0}(v)=1 \\[2pt]
    \ \varphi_{i}(v)=v_{i}, ~1 \leq i \leq 3 \\[2pt]
    \ \varphi_{4}(v)=\frac{1}{\sqrt{6}}(|v|^2 -3)
    \end{cases}
\end{equation}
with $\langle \varphi_{i}M,\varphi_{j}M \rangle_{L_{v}^{2}}=\delta_{ij}$. Let $h_{\varepsilon}$ be a solution to \eqref{linearized_Boltzmann_equation_Euler_scaling}, we define the fluid quantities associated to $h_{\varepsilon}$ by 
\begin{equation}\label{fluid quantities of h}
    \begin{cases}
    \ \rho_{\varepsilon}(t,x) =\langle h_{\varepsilon}, \varphi_{0}M \rangle_{L_{v}^{2}} \\[2pt]
    \ u_{\varepsilon}(t,x) =\langle h_{\varepsilon}, vM \rangle_{L_{v}^{2}} \\[2pt]
    \ T_{\varepsilon}(t,x) =\langle h_{\varepsilon}, \varphi_{4}M \rangle_{L_{v}^{2}}, 
\end{cases}
\end{equation}

\noindent which are precisely the coefficients of $\pi_{\mathcal{L}}(h_{\varepsilon})$ with respect to the orthonormal basis. We now describe the micro-macro decomposition for \eqref{linearized_Boltzmann_equation_Euler_scaling}. 

We first decompose $h_{\varepsilon}$ as 
\begin{equation}\label{MM_decomp_linearized_Boltzmann}
    h_{\varepsilon} = \pi_{\mathcal{L}}(h_{\varepsilon}) + \varepsilon g_{\varepsilon}=\bm{m}_{\varepsilon}^{T}\bm{\varphi}M + \varepsilon g_{\varepsilon}, 
\end{equation}
\noindent where $\bm{m}_{\varepsilon}=(\rho_{\varepsilon},u_{1,\varepsilon},u_{2,\varepsilon},u_{3,\varepsilon},T_{\varepsilon})^T , \bm{\varphi}=(\varphi_{0},\varphi_{1},\varphi_{2},\varphi_{3},\varphi_{4})^T$. To simplify the notation, we delete the subscript $\varepsilon$ in the expression above and let $\tilde{h}=\bm{m}^{T}\bm{\varphi}M$. Then the decomposition becomes $h=\tilde
{h}+\varepsilon g$. We substitute this ansatz into \eqref{linearized_Boltzmann_equation_Euler_scaling} to get 
\begin{equation*}
    \partial_{t}(\tilde{h}+\varepsilon g) + v \cdot \nabla_{x}(\tilde{h}+\varepsilon g)=\frac{1}{\varepsilon}\mathcal{L}(\tilde{h}+\varepsilon g)
\end{equation*} 

\noindent which is equivalent to 
\begin{equation}\label{MM_linearized_Boltzmann_1}
    \partial_{t}\tilde{h} + \varepsilon \partial_{t}g + v \cdot \nabla_{x}\tilde{h}+ \varepsilon v \cdot \nabla_{x}g = \mathcal{L}(g). 
\end{equation} 

Taking $\pi_{\mathcal{L}}$ on both sides of \eqref{MM_linearized_Boltzmann_1} gives the macroscopic part of the equation
\begin{equation}\label{MM_linearized_Boltzmann_macro}
    \partial_{t}\tilde{h}+\pi_{\mathcal{L}}(v \cdot \nabla_{x}\tilde{h})+ \varepsilon \pi_{\mathcal{L}}(v \cdot \nabla_{x}g)=0, 
\end{equation}

\noindent where we used $\pi_{\mathcal{L}}(g)=0$,  $\langle \mathcal{L}(f), \varphi_{i}M \rangle_{L_{v}^{2}}=0$,  \, $\forall f \in L^2_{v}, ~0 \leq i \leq 4$. Note that \eqref{MM_linearized_Boltzmann_macro} is equivalent to 
\begin{equation}\label{MM_linearized_Boltzmann_macro_2}
    \partial_{t}\bm{m}+\nabla_{x} \cdot \langle v\tilde{h}\bm{\varphi}M \rangle + \varepsilon \nabla_{x} \cdot \langle vg\bm{\varphi}M \rangle =0, 
\end{equation}

\noindent where $\langle f \rangle := \int_{\mathbb{R}^3} f(v) ~dv$ and $\langle v\tilde{h}\bm{\varphi}M \rangle$ stands for $(\langle v\tilde{h} \varphi_{0} M \rangle, \langle v\tilde{h} \varphi_{1} M \rangle,..., \langle v\tilde{h} \varphi_{4} M \rangle)^{T}$ (similarly for $\langle v g \bm{\varphi}M \rangle$). Taking $I-\pi_{\mathcal{L}}$ on both sides of \eqref{MM_linearized_Boltzmann_1} leads to the microscopic part 
\begin{equation}\label{MM_linearized_Boltzmann_micro}
    \varepsilon \partial_{t}g + (I-\pi_{\mathcal{L}})v \cdot \nabla_{x}\tilde{h}+ \varepsilon (I-\pi_{\mathcal{L}})v \cdot \nabla_{x}g = \mathcal{L}(g). 
\end{equation}

The coupled system \eqref{MM_linearized_Boltzmann_macro_2}-\eqref{MM_linearized_Boltzmann_micro} is the kinetic/fluid formulation of \eqref{linearized_Boltzmann_equation_Euler_scaling}, and we can easily recover \eqref{linearized_Boltzmann_equation_Euler_scaling} from the system by taking the dot product of \eqref{MM_linearized_Boltzmann_macro_2} with $\bm{\varphi}M$, and then add the result to \eqref{MM_linearized_Boltzmann_micro}. 

As $\varepsilon \rightarrow 0$, the coupled system becomes 

\begin{equation}\label{fluid limit of linearized Boltzmann}
    \begin{cases}
    \ \partial_{t}\bm{m}+\nabla_{x} \cdot \langle v\tilde{h}\bm{\varphi}M \rangle \rangle =0,  \\[2pt]
    \ (I-\pi_{\mathcal{L}})v \cdot \nabla_{x}\tilde{h} = \mathcal{L}(g). 
\end{cases}
\end{equation}

\noindent Thus $g=\mathcal{L}^{-1}((I-\pi_{\mathcal{L}})v \cdot \nabla_{x}\tilde{h})$ and $(\rho, u, T)$ satisfies $$\partial_{t}\bm{m}+\nabla_{x} \cdot \langle v\tilde{h}\bm{\varphi}M \rangle \rangle =0, $$
which takes the following explicit form 

\begin{equation}\label{acoustic equations}
    \begin{cases}
        \ \partial_{t}\rho + \nabla_{x} \cdot u = 0\\
        \ \partial_{t}u + \nabla_{x}(\rho +\frac{\sqrt{6}}{3}T)=0\\
        \ \partial_{t}T + \nabla_{x} \cdot (\frac{\sqrt{6}}{3}u) =0 
    \end{cases}
\end{equation}

\noindent This is precisely the acoustic equations as described in \cite{JJG2010}, with a slight change in the constant coefficients due to different choices of the orthonormal basis. Hence our macro-micro formulation correctly captures the fluid limit of \eqref{linearized_Boltzmann_equation_Euler_scaling}. 

Our goal is to obtain $(\rho,u,T,g)$ using neural networks. We will denote $f_{\theta}$ to be the approximation of $f$ using a neural network with parameters $\theta$. In the APNN framework, we will construct networks $(\rho_{\theta},u_{\theta},T_{\theta},g_{\theta})$ to approximate $(\rho,u,T,g)$, and we refer to \cite{LiuZhu2024} for a detailed description of the APNN structure. Notice that the function $g$ contains the velocity variable $v$, which has an unbounded domain $\mathbb{R}^3$. But when using neural networks to approximate functions, we typically consider functions in a bounded domain. So we will restrict our training of $g_{\theta}$ over some bounded domain $\mathcal{D} \subset \mathbb{R}^3$. For technical purpose, we then extend $g_{\theta}$ over the entire velocity space by setting $g_{\theta}=0$ for $v \in \mathbb{R}^3-\mathcal{D}$. Moreover, $\pi_{\mathcal{L}}(g)=0$ if $g$ is the real solution to \eqref{MM_linearized_Boltzmann_micro}. This clearly does not necessarily hold for an arbitrary $g_{\theta}$. Hence, we perform a post-processing of our network to replace $g_{\theta}$ by $g_{\theta} - \pi_{\mathcal{L}}(g_{\theta})$ to make sure $\pi_{\mathcal{L}}(g_{\theta})=0$ and hence $\pi_{\mathcal{L}}(h_{\theta}) = \tilde{h}_{\theta}$, where 
\begin{equation*}
    \begin{cases}
     \ \bm{m}_{\theta} = (\rho_{\theta},u_{\theta},T_{\theta})^{T} \\[2pt]
     \ \tilde{h}_{\theta} = \bm{m}_{\theta}^{T}\bm{\varphi}M \\[2pt]
     \ h_{\theta} = \tilde{h}_{\theta} + \varepsilon g_{\theta} 
    \end{cases}
\end{equation*}
The best approximation $(\rho_{\theta},u_{\theta},T_{\theta},g_{\theta})$ is obtained by minimizing the loss function of our neural network, which we now describe. Replacing $(\rho,u,T,g)$ by $(\rho_{\theta},u_{\theta},T_{\theta},g_{\theta})$ in \eqref{MM_linearized_Boltzmann_macro_2} and \eqref{MM_linearized_Boltzmann_micro}, we get 
\begin{equation}\label{MM_formulation_NN}
    \begin{cases}
        \ \partial_{t}\bm{m}_{\theta}+\nabla_{x} \cdot \langle v\tilde{h}_{\theta}\bm{\varphi}M \rangle + \varepsilon \nabla_{x} \cdot \langle vg_{\theta}\bm{\varphi}M \rangle= \bm{d}_{\theta}^{1} \\[2pt]
        \ \varepsilon \partial_{t}g_{\theta} + (I-\pi_{\mathcal{L}})v \cdot \nabla_{x}\tilde{h}_{\theta} + \varepsilon (I-\pi_{\mathcal{L}})v \cdot \nabla_{x}g_{\theta} - \mathcal{L}(g_{\theta}) =d_{\theta}^{2}
    \end{cases}
\end{equation}

We further define 
\begin{equation}\label{loss function initial}
    d_{\theta}^{ini} = h_{\theta}(0,x,v) - h(0,x,v)
\end{equation}
\begin{equation}\label{loss function boundary}
    \begin{aligned}
        \ d_{\theta}^{b} = &(h_{\theta}(t,x_1,x_2,\pi,v) - h_{\theta}(t,x_1,x_2,-\pi,v))^{2} + (h_{\theta}(t,x_1,\pi,x_3,v) \\
        \ &- h_{\theta}(t,x_1,-\pi,x_3,v))^{2} + (h_{\theta}(t,\pi,x_2,x_3,v) - h_{\theta}(t,-\pi,x_2,x_3,v))^{2}
    \end{aligned}
\end{equation}

\noindent which correspond to the losses on the initial data and boundary value, respectively. We define the generalization error by
\begin{equation}\label{Loss APNN}
    \mathcal{R}^{G}_{\theta} = \mathcal{R}_{\theta}^{1} +\mathcal{R}_{\theta}^{2} + \mathcal{R}_{\theta}^{ini}+ \mathcal{R}_{\theta}^{b}
\end{equation}

\noindent where 
\begin{equation}\label{loss composition}
    \begin{cases}
        \ \mathcal{R}_{\theta}^{1} = \int_{0}^{T} \mathcal{R}_{\theta,t}^{1} ~dt, ~~\mathcal{R}_{\theta,t}^{1} = \int_{\mathbb{T}^3} |\bm{d}^1_{\theta}|^{2} ~dx \\[2pt]
        \ \mathcal{R}_{\theta}^{2} = \int_{0}^{T} \mathcal{R}_{\theta,t}^{2} ~dt,~~\mathcal{R}_{\theta,t}^{2} = \int_{\mathbb{T}^3} \int_{\mathcal{D}} |d^2_{\theta}|^{2} ~dvdx\\[2pt]
        \ \mathcal{R}_{\theta}^{ini} = \int_{\mathbb{T}^3} \int_{\mathcal{D}} |d^{ini}_{\theta}|^{2} ~dvdx\\[2pt]
        \ \mathcal{R}_{\theta}^{b} = \int_{0}^{T} \mathcal{R}_{\theta,t}^{b} ~dt, ~~\mathcal{R}_{\theta,t}^{b} = \int_{\partial\mathbb{T}^3} \int_{\mathcal{D}} |d^{b}_{\theta}|^{2} ~dvd\sigma(x). 
    \end{cases}    
\end{equation}

\noindent Notice that as $\epsilon \rightarrow 0$, the residual loss \eqref{MM_formulation_NN} will tend to the residual loss of the fluid limit \eqref{fluid limit of linearized Boltzmann}, which clearly demonstrates the AP property of the network. 

Since we will eventually prove convergence results in the $H^1$ space, we need a $H^1$-counterpart of the loss function. Hence we further define $d^{1,\nabla_{x}}_{\theta} = \nabla_{x}d^{1}_{\theta}, ~d^{2,\nabla_{x}}_{\theta} = \nabla_{x}d^{2}_{\theta}$, and let  $d^{ini,\nabla_{x}}_{\theta}, d^{b,\nabla_{x}}_{\theta}$ similarly as $d_{\theta}^{ini}, d_{\theta}^{b}$, replacing $h_{\theta}, h$ by $\nabla_{x}h_{\theta}, \nabla_{x}h$. We then define $\mathcal{R}_{\theta}^{G,\nabla_{x}}$ similarly as $\mathcal{R}_{\theta}^{G}$, replacing each integrand by the corresponding $\nabla_x$-counterpart, and we define $\mathcal{R}_{\theta}^{G,\nabla_{v}}$ in the same way. Finally, we let
\begin{equation}\label{loss APNN H1}
    \mathcal{R}^{G,H^1}_{\theta} = \mathcal{R}^{G}_{\theta} + \mathcal{R}^{G,\nabla_{x}}_{\theta} + \mathcal{R}^{G, \nabla_{v}}_{\theta}, 
\end{equation}

\noindent which is the final loss function we are going to minimize during the training process of our networks.

\subsection{APNN for the SG system}

We now generalize the previous APNN framework to linearized Boltzmann equation with uncertainty. We assume again acoustic scaling and then the equation reads 
\begin{equation}\label{fully linear Botzmann with UQ}
    \begin{cases}
    \ \partial_{t}h + v\cdot \nabla_{x}h = \frac{1}{\varepsilon}\mathcal{L}(h),  \\[2pt]
    \ h(0,x,v,z)= h_{I}. 
\end{cases}
\end{equation}

To deal with the random variale $z$, we follow the stochastic Galerkin method introduced in \ref{SG method}. Let $h \approx h_K =\sum_{i=1}^{K}h^{i}\phi_{i}(z)$, and we again write $h=\tilde{h}+\varepsilon g=\bm{m}^{T}\bm{\varphi}M+\varepsilon g$. Then 
\begin{equation*}
    \begin{aligned}
        \ h^{i}(t,x,v)=& \int_{I_z} h\phi_{i}(z)\pi(z) ~dz \\
        \ =& \int_{I_z} (\bm{m}^{T}\bm{\varphi}M+\varepsilon g)\phi_{i}(z)\pi(z) ~dz \\
        \ =& \int_{I_z} \bm{m}^{T}\bm{\varphi}M\phi_{i}(z)\pi(z) ~dz + \varepsilon \int_{I_z} g\phi_{i}(z)\pi(z) ~dz \\
        \ =& (\bm{m}^{i})^{T}\bm{\varphi}M + \varepsilon g^{i}, 
    \end{aligned}
\end{equation*}
\noindent where $\bm{m}^{i}=(\rho^{i},u^{i}, T^{i})^{T}$, which can be easily proved to be the coefficients of $\pi_{\mathcal{L}}(h^i)$ with respect to the orthonormal basis $\{ \bm{\varphi}M \}$, due to the fact that $\pi_{\mathcal{L}}$ is interchangeable with $\int_{I_z} \cdot ~dz$. Hence we also have $g^{i} \in N(\mathcal{L})^{\perp}$ for $1 \leq i \leq K$. 

As in \eqref{SG_linearized_Boltzmann_equation}, the stochastic Galerkin system for \eqref{fully linear Botzmann with UQ} is given by 
\begin{equation}\label{SG system for fully linear Boltzmann}
    \partial_{t}h^i + v\cdot \nabla_{x}h^i = \frac{1}{\varepsilon}\mathcal{L}_{i}(h_K), 
\end{equation}
\noindent for $1 \leq i \leq K$ where $\mathcal{L}_{i}(h_K)$ is defined in \eqref{SG_linear_term}. Substituting $h^i = (\bm{m}^{i})^{T}\bm{\varphi}M + \varepsilon g^{i}$ into \eqref{SG system for fully linear Boltzmann} and taking $\pi_{\mathcal{L}}$ and $I - \pi_{\mathcal{L}}$ respectively on both sides of the equation will lead to 
\begin{equation}\label{MM_formulation_SG}
    \begin{cases}
        \ \partial_{t}\bm{m}^{i} + \nabla_{x} \cdot \langle v\tilde{h}^{i}\bm{\varphi}M \rangle + \varepsilon \nabla_{x} \cdot \langle vg^{i}\bm{\varphi}M \rangle = 0,   \\[2pt]
        \ \varepsilon \partial_{t}g^{i} + (I-\pi_{\mathcal{L}})\nabla_{x} v \cdot \tilde{h}^{i} + \varepsilon (I - \pi_{\mathcal{L}}) v \cdot \nabla_{x}g^{i} = \mathcal{L}_{i}(g_K), 
    \end{cases}
\end{equation}
\noindent for $1 \leq i \leq K$. This is the macro-micro formulation for \eqref{SG system for fully linear Boltzmann}. 

As in the deterministic case, we will construct networks $(\rho_{\theta}^{i},u_{\theta}^{i},T_{\theta}^{i},g_{\theta}^{i})$ to approximate $(\rho_{\theta}^{i},u_{\theta}^{i},T_{\theta}^{i},g_{\theta}^{i})$. The losses $(\bm{d}_{\theta}^{1,i},d_{\theta}^{2,i},d_{\theta}^{i,i},d_{\theta}^{b,i})$ are defined similarly as $(\bm{d}_{\theta}^{1},d_{\theta}^{2},d_{\theta}^{i},d_{\theta}^{b})$, replacing $(\rho_{\theta},u_{\theta},T_{\theta},g_{\theta})$ by $(\rho_{\theta}^{i},u_{\theta}^{i},T_{\theta}^{i},g_{\theta}^{i})$. Then 
\begin{equation}\label{MM_formulation_SG_NN}
    \begin{cases}
        \ \partial_{t}\bm{m}^{i}_{\theta} + \nabla_{x} \cdot \langle v\tilde{h}^{i}_{\theta}\bm{\varphi}M \rangle + \varepsilon \nabla_{x} \cdot \langle vg^{i}_{\theta}\bm{\varphi}M \rangle = \bm{d}_{\theta}^{1,i},   \\[2pt]
        \ \varepsilon \partial_{t}g^{i}_{\theta} + (I-\pi_{\mathcal{L}}) v \cdot \nabla_{x}\tilde{h}^{i}_{\theta} + \varepsilon (I - \pi_{\mathcal{L}}) v \cdot \nabla_{x}g^{i}_{\theta} - \mathcal{L}_{i}(g_{K,\theta})= d_{\theta}^{2,i}. 
    \end{cases}
\end{equation}
The generalization error in the stochastic case is defined by 
\begin{equation}\label{APNN loss stochastic}
    \mathcal{R}^{G}_{\theta,sto} = \mathcal{R}_{\theta,sto}^{1} +\mathcal{R}_{\theta,sto}^{2} + \mathcal{R}_{\theta,sto}^{i}+ \mathcal{R}_{\theta,sto}^{b}, 
\end{equation}

\noindent where 
\begin{equation*}
    \mathcal{R}_{\theta,sto}^{1} = \sum_{i=1}^{K}i^{2q}\int_{0}^{T} \mathcal{R}_{\theta,t}^{1,i} ~dt, ~~\mathcal{R}_{\theta,t}^{1,i} = \int_{\mathbb{T}^3} |\bm{d}^{1,i}_{\theta}|^{2} ~dx, 
\end{equation*}
\noindent and $\mathcal{R}_{\theta,sto}^{2}, \mathcal{R}_{\theta,sto}^{i},  \mathcal{R}_{\theta,sto}^{b}$ are defined similarly. Finally, we define 
\begin{equation}\label{loss APNN H1 stochastic}
    \mathcal{R}^{G,H^1}_{\theta,sto} = \mathcal{R}^{G}_{\theta,sto} + \mathcal{R}^{G,\nabla_{x}}_{\theta,sto} + \mathcal{R}^{G, \nabla_{v}}_{\theta,sto}\,, 
\end{equation}
\noindent where each term is defined analogously to the deterministic case. 
\section{Main results}

\subsection{Existence of APNN with arbitrarily small loss}
Our first result states that there exists a neural network so that the APNN loss $\mathcal{R}_{\theta,sto}^{G,H^1}$ is sufficiently small. Notice that $g^{i}_{\theta}$ can only approximate $g^i$ in a bounded subset of the velocity space. Hence inspired by \cite{Yangxu2024},we need to make some technical assumptions on $g^i$ outside some bounded subset of $\mathbb{R}^3$ so that the approximation is valid over the entire velocity space. To be specific, we define the following quantities.  

Let 
\begin{equation*}
    c_{i} = \int_{0}^{T}\int_{\mathbb{T}^3}\int_{\mathbb{R}^3-\mathcal{D}} |g^{i}|^2 + |v\cdot \nabla_{x}g^{i}|^2 ~dvdxdt 
\end{equation*}

\noindent Define $c_{i}^{\nabla_{x}}, c_{i}^{\nabla_{v}}$ in a similarly way, replacing $g^i$ by $\nabla_{x}g^i$ and $\nabla_{v}g^i$ respectively in the definition of $c_{i}$. Let 
\begin{equation*}
    c_{ik}^{\Lambda} = \int_{0}^{T}\int_{\mathbb{T}^3}\int_{\mathbb{R}^3-\mathcal{D}} |\nu_{ik}\nabla_{v}g^{k}|^2 + |(\nabla_{v}\nu_{ik})g^{k}|^2 ~dvdxdt. 
\end{equation*}

\noindent where $\nu_{ik}$ is defined similarly as \eqref{collision_frequency}, replacing the collision kernel with the proper kernels in the SG setting. 

We make the following assumption, which asserts that quantities related to the microscopic part $\bm{g}$ outside some sufficiently large $\mathcal{D}$ is negligible. 

\bigskip
\noindent \textbf{Assumption 1} For any $\eta >0$, there exists a bounded domain $\mathcal{D} \subset \mathbb{R}^3$ large enough, such that $c_{\mathcal{D}}:= \sum_{i=1}^{K}(c_{i} + c_{i}^{\nabla_{x}} + c_{i}^{\nabla_{v}} + \sum_{k=1}^{K}c_{ik}^{\Lambda}) < \eta$. 
\begin{theorem}
    Consider the linearized Boltzmann equation with uncertainty \eqref{fully linear Botzmann with UQ}, and let $(\bm{\rho},\bm{u},\bm{\theta},\bm{g})$ be the solution to the macro-micro formulation \eqref{MM_formulation_SG} such that the microscopic part $\bm{g}$ satisfies Assumption 1. Then for any $\delta>0$, there exists a bounded domain $\mathcal{D} \subset \mathbb{R}^3$ and a network $(\bm{\rho}_{\theta},\bm{u}_{\theta},\bm{T}_{\theta},\bm{g}_{\theta})$ such that 
    $$\mathcal{R}_{\theta,sto}^{G,H^1} < \delta. $$
\end{theorem}

\textit{Proof:} To prove the statement, we need the following lemma from \cite{LiuZhu2024}, which is basically the Universal Approximation Theorem for neural networks. 

\begin{lemma}
    Let $\Omega \subset \mathbb{R}^N$ be a bounded subset. Suppose $f \in C^{2}(\Omega)$. Then for any $\eta >0$, there exists a two-layer neural network $f_{\theta}$ such that, 

    \begin{equation*}
        ||f-f_{\theta}||_{W^{2,\infty}(\Omega)} < \eta 
    \end{equation*}
\end{lemma}

Let $\mathcal{D}$ be a bounded domain in $\mathbb{R}^3$ such that $c_{\mathcal{D}}<d\delta$ for some $d>0$ to be chosen later. Let $\Omega=[0,T] \times \mathbb{T}^3 \times \mathcal{D}$. Suppose $h^{i} \in C^2(\Omega)$ for all $1 \leq i \leq K$. Then from Lemma 1, for any $\eta>0$, there exists a network $(\bm{\rho}_{\theta},\bm{u}_{\theta},\bm{T}_{\theta},\bm{g}_{\theta})$
such that 
\begin{equation}\label{NN_approximation error}
    \begin{cases}
        \ ||\rho^{i}_{\theta} - \rho^{i}||_{W^{2,\infty}(\Omega)} < \eta, ~||u^{i}_{\theta} - u^{i}||_{W^{2,\infty}(\Omega)} < \eta, \\[2pt]
        \ ||T^{i}_{\theta} - T^{i}||_{W^{2,\infty}(\Omega)} < \eta, ~||g^{i}_{\theta} - g^{i}||_{W^{2,\infty}(\Omega)} < \eta. 
    \end{cases}
\end{equation}

Combining \eqref{MM_formulation_SG} and \eqref{MM_formulation_SG_NN}, we have 
\begin{equation}\label{MM_formulation_SG_NN_rewrite}
    \begin{cases}
        \ \partial_{t}\hat{\bm{m}}^{i}_{\theta} + \nabla_{x} \cdot \langle v\hat{\tilde{h}}^{i}_{\theta}\bm{\varphi}M \rangle + \varepsilon \nabla_{x} \cdot \langle v\hat{g}^{i}_{\theta}\bm{\varphi}M \rangle = -\bm{d}_{\theta}^{1,i}\\
        \varepsilon \partial_{t}\hat{g}^{i}_{\theta} + (I-\pi_{\mathcal{L}}) v \cdot \nabla_{x}\hat{\tilde{h}}^{i}_{\theta} + \varepsilon (I - \pi_{\mathcal{L}}) v \cdot \nabla_{x}\hat{g}^{i}_{\theta} - \mathcal{L}_{i}(\hat{g}_{K,\theta})=- d_{\theta}^{2,i}
    \end{cases}
\end{equation}

For the rest of the proof, we use the notation $||\cdot||=||\cdot||_{L^2_{[0,T]\times \mathbb{T}^3\times\mathbb{R}^3}}$, $||\cdot||_{\Omega} = ||\cdot||_{L^2_{\Omega}}$ and $||\cdot||_{\Omega-\mathcal{D}} = ||\cdot||_{L^2_{[0,T]\times \mathbb{T}^3\times (\mathbb{R}^3 - \mathcal{D})}}$. Consider $\mathcal{R}_{\theta}^{2,i} = ||d_{\theta}^{2,i}||_{\Omega}^{2}$. From \eqref{NN_approximation error}, 
\begin{equation*}
    \varepsilon ||\partial_{t}\hat{g}^{i}_{\theta}||_{\Omega}^{2} < ||\partial_{t}\hat{g}^{i}_{\theta}||_{\Omega}^{2} \leq |\Omega|\cdot||\partial_{t}\hat{g}^{i}_{\theta}||_{L_{\Omega}^{\infty}}^{2} < |\Omega|\eta^2
\end{equation*}
 and 
\begin{equation*}
    \begin{aligned}
        \ ||\varepsilon (I - \pi_{\mathcal{L}}) v \cdot \nabla_{x}\hat{g}^{i}_{\theta}||_{\Omega}^{2} \leq ||(I - \pi_{\mathcal{L}}) v \cdot \nabla_{x}\hat{g}^{i}_{\theta}||^2 \leq& ||v \cdot \nabla_{x}\hat{g}^{i}_{\theta}||^2 \\
        \ =& ||v \cdot \nabla_{x}\hat{g}^{i}_{\theta}||^2_{\Omega} + ||v \cdot \nabla_{x}\hat{g}^{i}_{\theta}||^2_{\Omega-\mathcal{D}} \\
        \ <& C_1\eta^2 + ||v \cdot \nabla_{x}g^{i}||^2_{\Omega-\mathcal{D}} \\
        \ \leq& C_1\eta^2 + d\delta
    \end{aligned}
\end{equation*}

\noindent where we use the fact that $g^{i}_{\theta}$ vanishes when $v \in \mathbb{R}^3 - \mathcal{D}$ and assumption 1 in the last two inequalities. Similarly, we can show $||(I-\pi_{\mathcal{L}}) v \cdot \nabla_{x}\hat{\tilde{h}}^{i}_{\theta}||^2_{\Omega} < C_2\eta^2 + d\delta$.

Finally, since $\mathcal{L}_{ik}$ is a bounded operator on $L^2(\mathbb{R}^3)$, for all $1 \leq i,k \leq K$, we have 
\begin{equation*}
    ||\mathcal{L}_{ik}(\hat{g}_{\theta}^{k})||^2_{\Omega} \leq ||L_{ik}(\hat{g}_{\theta}^{k})||^2 \leq C_3||\hat{g}_{\theta}^{k}||^2 = C_3(||\hat{g}_{\theta}^{k}||^2_{\Omega}+||{g}^{k}||^2_{\Omega-\mathcal{D}}) < C_3\eta^2 + C_3 d\delta
\end{equation*}

Collecting all the above inequalities, we get $\mathcal{R}_{\theta}^{2,i} < C\eta^2 + C^{'}d\delta$ for some $C, C^{'}>0$. Using the same arguments, we can show $\mathcal{R}_{\theta}^{2,i,\nabla_{x}} < C\eta^2 + C^{'}d\delta$. For $\mathcal{R}_{\theta}^{2,i,\nabla_{v}}$, we will apply similar arguments, with additional observations that 
\begin{equation*}
    \begin{aligned}
        \ ||\nabla_{v} (I - \pi_{\mathcal{L}}) v \cdot \nabla_{x}\hat{g}^{i}_{\theta}||^2 \leq& ||\nabla_{x}\hat{g}^{i}_{\theta}||^2 + ||v\cdot \nabla_{v}\nabla_{x}\hat{g}^{i}_{\theta}||^2 + ||\nabla_{v}\pi_{\mathcal{L}}(v \cdot \nabla_{x}\hat{g}^{i}_{\theta})||^2 \\
        \ \leq& ||\nabla_{x}\hat{g}^{i}_{\theta}||^2 + ||v\cdot \nabla_{v}\nabla_{x}\hat{g}^{i}_{\theta}||^2 + C_{\pi 1}||v \cdot \nabla_{x}\hat{g}^{i}_{\theta}||^2
    \end{aligned}
\end{equation*}

\noindent where we use \eqref{H_3_1} for the last inequality, and 
\begin{equation*}
    \begin{aligned}
        \ ||\nabla_{v}\mathcal{L}_{ik}(\hat{g}_{\theta}^{k})||^2 \leq & 2||\nabla_{v}K_{ik}(\hat{g}_{\theta}^{k})||^2 + 2||\nabla_{v}\Lambda_{ik}(\hat{g}_{\theta}^{k})||^2 \\
        \ \leq& 2\eta||\nabla_{v}\hat{g}_{\theta}^{k}||^2 +2C(\eta)||\hat{g}_{\theta}^{k}||^2 +4||(\nabla_{v}\nu_{ik})\hat{g}_{\theta}^{k}||^2 +4||\nu_{ik}\nabla_{v}\hat{g}_{\theta}^{k}||^2
    \end{aligned}
\end{equation*}

\noindent where we use \eqref{H_2} and the definition for $\Lambda_{ik}$. So we have $\mathcal{R}_{\theta,sto}^{2}+\mathcal{R}_{\theta,sto}^{2,\nabla_{x}}+\mathcal{R}_{\theta,sto}^{2,\nabla_{v}} < C\eta^2 + C^{'}d\delta$. The error estimates for $\mathcal{R}^{1}_{\theta,sto}, \mathcal{R}^{ini}_{\theta,sto}, \mathcal{R}^{b}_{\theta,sto}$ and their $\nabla_{x}(or ~\nabla_v)$-counterparts are analogous and simpler. So we omit the details here. Finally, setting $d$ and $\eta$ to be small enough will yield the desired result. 

\subsection{Convergence of APNN solution}

Our next theorem suggests that when the APNN loss $\mathcal{R}_{\theta,sto}^{G,H^1}$ is small enough, $(\bm{\rho}_{\theta},\bm{u}_{\theta},\bm{\theta}_{\theta},\bm{g}_{\theta})$ will tend towards the true solution $(\bm{\rho},\bm{u},\bm{\theta},\bm{g})$ in the weighted $H^1$-norm. Again, we need to make some technical assumptions on the microscopic part $\bm{g}$, which is given as follows. 

Let 
\begin{equation}
    \begin{cases}
        \ \tilde{R}_{1,t}^{i}=\int_{\mathbb{T}^3} \int_{\mathbb{R}^3 - \mathcal{D}}|\partial_{t}g^{i}|^2 +|v \cdot \nabla_{x}g^{i}|^2 + |\mathcal{L}_{i}(g_K)|^2 ~dvdx, ~\tilde{R}_{1}^{i} = \int_{0}^{T} \tilde{R}_{1,t}^{i} ~dt,  \\[2pt]
        \ \tilde{R}_{2}^{i} = \int_{\mathbb{T}^3} \int_{\mathbb{R}^3-\mathcal{D}} |g^{i}(0,x,v)|^{2} ~dvdx, \\[2pt]
        \ \tilde{R}_{3,t}^{i} = \int_{\partial\mathbb{T}^3} \int_{\mathbb{R}^3-\mathcal{D}} v|g^{i}|^{2} ~dvd\sigma(x), ~\tilde{R}_{3}^{i} = \int_{0}^{T} \tilde{R}_{3,t}^{i} ~dt. 
    \end{cases}
\end{equation}

Replacing $g$ by $\nabla_{x}g$ and $\nabla_{v}g$ respectively, we can define $\tilde{R}_{j}^{\nabla_x,i}$ and $\tilde{R}_{j}^{\nabla_v,i}$ in a similar way for $j=1,2,3$. We let $\tilde{R} = \sum_{i=1}^{K}\sum_{j=1}^{3} (\tilde{R}_{j}^{i}+ \tilde{R}_{j}^{\nabla_x,i}+\tilde{R}_{j}^{\nabla_v,i})$. We make the following assumption on $\tilde{R}$, 

\bigskip
\noindent \textbf{Assumption 2} For any $0 < \varepsilon <1$, there exists $\mathcal{D} \subset \mathbb{R}^3$ large enough such that $\tilde{R} < \delta(\varepsilon)$ with $\delta(\varepsilon) \rightarrow 0$ as $\varepsilon \rightarrow 0$. 
\begin{theorem}
    Let $(\bm{\rho}_{\theta},\bm{u}_{\theta},\bm{\theta}_{\theta},\bm{g}_{\theta})$ be approximations to $(\bm{\rho},\bm{u},\bm{T},\bm{g})$ using a neural network with parameters $\theta$ such that $\bm{g}$ satisfies assumption 2, and let $\hat{\bm{h}}_{\theta} = \bm{h} - \bm{h}_{\theta}$. Then for any $t \in (0,T]$ and $0 < \varepsilon < 1$, 
 \begin{equation*}
        E^{K}_{t}(\hat{\bm{h}}_{\theta}) \leq \frac{\tilde{C}(\mathcal{R}_{\theta,sto}^{G,H^1}+\delta(\varepsilon))}{\varepsilon} e^{-\varepsilon \tau t}, 
    \end{equation*}
\noindent for some $\tilde{C}, \tau >0$ independent of $\varepsilon$. 
\end{theorem}

\textit{Proof}. Combing the two equations in \eqref{MM_formulation_SG_NN_rewrite} gives 
\begin{equation}\label{Boltzmann_eq_SG_NN}
    \partial_{t}\hat{h}_{\theta}^{i} + v \cdot \nabla_{x}\hat{h}_{\theta}^{i} = \frac{1}{\varepsilon}\mathcal{L}_{i}(\hat{h}_{K,\theta}) - \bm{d}^{1,i}_{\theta} \cdot \bm{\varphi}M - d^{2,i}_{\theta}
\end{equation}
\noindent where $\hat{h}^{i}_{\theta} = h^{i} - h^{i}_{\theta} ~~\forall 1 \leq i \leq K$. We let $A^{i}_{\theta} = - \bm{d}^{1,i}_{\theta} \cdot \bm{\varphi}M - d^{2,i}_{\theta}$. Unless specified otherwise, we use the notation $||\cdot||^2=||\cdot||^2_{L^2_{x,v}}$ and $\langle \cdot,\cdot \rangle = \langle \cdot,\cdot \rangle_{L^2_{x,v}}$ for the rest of the proof. Define 
\begin{equation*}
    \begin{cases}
        \ r_{1,t}^{i} := \int_{\mathbb{T}^3} \int_{\mathbb{R}^3 - \mathcal{D}}|d_{\theta}^{2,i}|^2 dvdx \\[2pt]
        \ r_{2}^{i} := \int_{\mathbb{T}^3} \int_{\mathbb{R}^3-\mathcal{D}} |h^{i}(0,x,v)|^{2} ~dvdx\\[2pt]
        \ r_{3,t}^{i} :=\int_{\partial\mathbb{T}^3} \int_{\mathbb{R}^3-\mathcal{D}} v|h^{i}|^{2} ~dvd\sigma(x). 
    \end{cases}
\end{equation*}
\noindent Since $g_{\theta}$ vanishes on $\mathbb{R}^3-\mathcal{D}$ and the Maxwellian $\mathcal{M}(v)$ is sufficiently small outside some large $\mathcal{D}$, we can derive that $\sum_{i=1}^{K}\int_{0}^{T} r_{1,t}^{i} ~dt < \delta(\varepsilon)$ using assumption 2. Similarly, $ \sum_{i=1}^{K}r_{2}^{i} < \delta(\varepsilon)$ and $\sum_{i=1}^{K}\int_{0}^{T} r_{3,t}^{i} ~dt < \delta(\varepsilon)$. Same arguments also hold for $r_{j}^{\nabla_{x},i}$ and $r_{j}^{\nabla_{v},i}$ with $1 \leq j \leq 3$. We then follow the idea in \cite{briant2015} to construct a Lyapunov functional $||\cdot||^2_{\varepsilon \perp}$ on $H^{1}_{x,v}$ that is equivalent to the standard Sobolev norm by 
\begin{equation}\label{def of Lyapunov}
    ||h||^2_{\varepsilon \perp} = a_1||h||^2 + a_2||\nabla_{x}h||^2 + a_3||\nabla_{v}h^{\perp}||^2 + a_4\varepsilon \langle  \nabla_{x}h, \nabla_{v}h\rangle ~~\forall h \in H^{1}_{x,v}
\end{equation}
\noindent for $a_1,a_2,a_3,a_4>0$. 

\begin{lemma}
    For any fixed choices of $a_1,a_3,a_4>0$ and $0 < \varepsilon < 1$, we can choose $a_2$ to be big enough so that $||\cdot||^2_{\varepsilon \perp}$ is equivalent to the standard Sobolev norm $||\cdot||^{2}_{H^{1}_{x,v}}$, with equivalence independent of $\varepsilon$. 
\end{lemma}

\textit{Proof of lemma.} Using Cauchy-Schwarz, Young's inequality and the fact that $\varepsilon <1 $, we have 
\begin{equation*}
    \begin{aligned}
         \ & a_4\varepsilon \langle  \nabla_{x}h, \nabla_{v}h\rangle \leq a_4\varepsilon\eta ||\nabla_{x}h||^2 + \frac{a_4\varepsilon}{\eta}||\nabla_{v}h||^2 \leq a_4\eta ||\nabla_{x}h||^2 + \frac{a_4}{\eta}||\nabla_{v}h||^2\\
         \ & a_4\varepsilon \langle  \nabla_{x}h, \nabla_{v}h\rangle \geq -a_4\varepsilon\eta ||\nabla_{x}h||^2 - \frac{a_4\varepsilon}{\eta}||\nabla_{v}h||^2 \geq -a_4\eta ||\nabla_{x}h||^2 - \frac{a_4}{\eta}||\nabla_{v}h||^2. 
    \end{aligned}
\end{equation*}

\noindent Write $h=h^{\perp}+\pi_{\mathcal{L}}(h)$ and using \eqref{H_3_1}, we get  
\begin{equation*}
    ||\nabla_{v}h||^2 \leq ||\nabla_{v}h^{\perp}||^2 + C_{\pi 1}||h||^2, 
\end{equation*}

\noindent so 
\begin{equation*}
    ||h||^2_{\varepsilon \perp} \geq (a_1-\frac{a_4C_{\pi 1}}{\eta})||h||^2 + (a_2-a_4\eta)||\nabla_{x}h||^2 + (a_3-\frac{a_4}{\eta})||\nabla_{v}h^{\perp}||^2. 
\end{equation*}
 
\noindent Thus, for any fixed choices $a_1,a_3,a_4 >0$, we can choose $\eta$ and $a_2$ big enough, so that $||h||^2_{\varepsilon \perp}>0 ~\forall h\neq 0$ and $||h||^2_{\varepsilon \perp} \geq C_1||h||^2_{\sim}$ for some $C_1>0$ where $||h||^2_{\sim} = ||h||^2+||\nabla_{x}h||^2+||\nabla_{v}h^{\perp}||^2$. Similarly, it can be proved that $||h||^2_{\varepsilon \perp} \leq C_2||h||^2_{\sim}$ for some $C_2>0$, independent of $\varepsilon$. So for any fixed $a_1,a_3,a_4 >0$ and $0 < \varepsilon < 1$, we can choose $a_2$ to be big enough such that $||\cdot||^2_{\varepsilon \perp}$ is equivalent to $||\cdot||^2_{\sim}$, with equivalence indenpendent of $\varepsilon$. 

Using \eqref{H_3_1}, we can derive 
\begin{equation*}
    ||\nabla_{v}h||^2 \leq ||\nabla_{v}h^{\perp}||^2 + C_{\pi 1}||h||^2 \leq ||\nabla_{v}h||^2 +2 C_{\pi 1}||h||^2. 
\end{equation*}

\noindent Hence $||\cdot||^2_{\sim}$ is equivalent to the standard Sobelov norm $||\cdot||^{2}_{H^{1}_{x,v}}$, and the result follows. 

Similarly, we define 
\begin{equation}\label{def of Lyapunov SG}
    E_{\varepsilon \perp,t}^{K}(\bm{h}) = a_1\sum_{i=1}^{K}||i^{q}h^i||^2 + a_2\sum_{i=1}^{K}||i^{q}\nabla_{x}h^i||^2 + a_3\sum_{i=1}^{K}||i^{q}\nabla_{v}h^{i,\perp}||^2 + a_4\varepsilon\sum_{i=1}^{K}i^{2q} \langle  \nabla_{x}h^i, \nabla_{v}h^i\rangle. 
\end{equation}

\noindent Then $E_{\varepsilon \perp,t}^{K}$ is equivalent to $E^{K}_{t}$ for proper choices of $a_1,a_2,a_3,a_4>0$. Next, we use four lemmas to bound the time evolution for each term in \eqref{def of Lyapunov SG}. 

\bigskip

\begin{lemma}
We have the estimate
    \begin{equation*}
        \partial_{t}\sum_{i=1}^{K}||i^{q}\hat{h}^{i}_{\theta}||^2 \leq -\frac{\lambda_{D}}{\varepsilon}\sum_{i=1}^{K} ||i^{q}\hat{h}^{i,\perp}_{\theta}||^{2}_{\Lambda} + R_1
    \end{equation*}

\noindent where 
\begin{equation*}
    R_1 = f_1\sum_{i=1}^{K}||i^qA^{i}_{\theta}||^2 + \frac{1}{f_1}\sum_{i=1}^{K}||i^q\hat{h}^{i}_{\theta}||^2 + C\sum_{i=1}^{K}i^{2q}\mathcal{R}_{\theta,t}^{b,i}+\sum_{i=1}^{K}i^{2q}r_{3,t}^{i} 
\end{equation*}

\noindent for $\lambda_{D}, f_1, C>0$, all independent of $\varepsilon$ with $f_1$ a free variable to be chosen later. 
\end{lemma}

\textit{Proof of lemma.} Take the inner product with respect to $\hat{h}^{i}_{\theta}$ on both sides of \eqref{Boltzmann_eq_SG_NN}, we get 
\begin{equation}\label{term_1}
    \begin{aligned}
        \ \frac{1}{2}\partial_{t}||\hat{h}^{i}_{\theta}||^2 = & -\langle v \cdot \nabla_{x}\hat{h}^{i}_{\theta}, \hat{h}^{i}_{\theta}\rangle + \frac{1}{\varepsilon}\sum_{k=1}^{K}\langle \mathcal{L}_{ik}(\hat{h}^{k}_{\theta}), \hat{h}^{i}_{\theta} \rangle + \langle A_{\theta}^{i}, \hat{h}^{i}_{\theta} \rangle. 
    \end{aligned}
\end{equation}

Notice that since $\hat{h}^{i}_{\theta}$ is not necessarily periodic, the operator $v\cdot \nabla_{x}$ is not skew-symmetric. Using integration by part and splitting $\int_{\mathbb{R}^3} = \int_{\mathcal{D}} + \int_{\mathbb{R}^3-\mathcal{D}}$, we have 
\begin{equation*}
    -\langle v \cdot \nabla_{x}\hat{h}^{i}_{\theta}, \hat{h}^{i}_{\theta}\rangle \leq C\mathcal{R}_{\theta,t}^{b,i}+r_{3,t}^{i}. 
\end{equation*}

To bound the term $\frac{1}{\varepsilon}\sum_{k=1}^{K}\langle \mathcal{L}_{ik}(\hat{h}^{k}_{\theta}), \hat{h}^{i}_{\theta} \rangle$, we follow exactly the same arguments in \cite{LiuDaus2019}, which we summarize as follows. Define 
\begin{equation}\label{Term-I}
    \textmd{Term I} = \sum_{i=1}^{K}i^{2q}\sum_{k=1}^{K}\langle \mathcal{L}_{ik}(\hat{h}^{k}_{\theta}), \hat{h}^{i}_{\theta} \rangle. 
\end{equation}

Let $\Theta_{i}=\frac{\hat{h}^{i,'}_{\theta,*}}{M^{'}_{*}}+\frac{\hat{h}^{i,'}_{\theta}}{M^{'}}-\frac{\hat{h}^{i}_{\theta,*}}{M_{*}}-\frac{\hat{h}^{i}_{\theta}}{M}$, $\Tilde{\Theta}_{i}=i^{q}\Theta_{i}$. Consider the change of variables $(v,v_{*}) \rightarrow (v^{'},v^{'}_{*}),(v_{*},v)$ respectively. It can be shown that 
\begin{equation*}
    \begin{aligned}
        \ \textmd{Term I}=&\sum_{i,k=1}^{K}i^{2q}\langle \mathcal{L}_{ik}(\hat{h}^{k}_{\theta}), \hat{h}^{i}_{\theta} \rangle \\
        \ =& -\frac{1}{4}\int_{\mathbb{T}^3}\int_{\mathbb{R}^3}\int_{\mathbb{R}^3 \times \mathbb{S}^2}  \mathcal{M}\mathcal{M}_{*} \times \textmd{Term A} ~~dv_{*}d\sigma dvdx, 
    \end{aligned}
\end{equation*}

\noindent where 
\begin{equation}\label{Term-A}
    \textmd{Term A}= \sum_{i,k=1}^{K}(\frac{i}{k})^{q}S_{ik}\Tilde{\Theta}_{i}\Tilde{\Theta}_{k}. 
\end{equation}

\noindent Using assumptions \eqref{assumptions on collision kernel}, \eqref{b linear in z} and \eqref{assumption on b0} on the collision kernel, it can be proved that 

\begin{equation*}
    \begin{aligned}
        \ \textmd{Term A} \geq \left(b_0 - (2^q+2)|b_1|C_z\right)\sum_{i=1}^{K}\Tilde{\Theta_{i}^{2}}  \geq D(\cos{\theta})\sum_{i=1}^{K}\Tilde{\Theta_{i}^{2}}, 
    \end{aligned}
\end{equation*}

\noindent and hence 

\begin{equation*}
    \begin{aligned}
        \ \textmd{Term I} \leq & -\frac{1}{4}\sum_{i=1}^{K}\int_{\mathbb{T}^3}\int_{\mathbb{R}^3}\int_{\mathbb{R}^3 \times \mathbb{S}^2}  \mathcal{M}\mathcal{M}_{*}\phi(|v-v_{*}|)D(\cos{\theta})\Tilde{\Theta_{i}^{2}} ~~dv_{*}d\sigma dvdx \\
        \ = & \sum_{i=1}^{K} i^{2q} \langle \mathcal{L}^{D}(\hat{h}^{i}_{\theta}),\hat{h}^{i}_{\theta} \rangle 
         \leq  -\lambda_{D}\sum_{i=1}^{K} i^{2q}||\hat{h}^{i,\perp}_{\theta}||^{2}_{\Lambda}, 
    \end{aligned}
\end{equation*}

\noindent where 

\begin{equation*}
    \mathcal{L}^{D}(h) = M\int_{\mathbb{R}^d \times \mathbb{S}^{d-1}} \phi(|v-v_{*}|)D(\cos{\theta})\mathcal{M}_{*} \left(\frac{h^{'}_{*}}{M^{'}_{*}} + \frac{h^{'}}{M^{'}} - \frac{h_{*}}{M_{*}} - \frac{h}{M}\right) dv_{*}d\sigma. 
\end{equation*}

\noindent We use the Cauchy-Schwarz and Young's inequality to bound $\langle A^{i}_{\theta}, \hat{h}^{i}_{\theta} \rangle$ by
\begin{equation*}
         \langle A^{i}_{\theta}, \hat{h}^{i}_{\theta} \rangle \leq f_1||A^{i}_{\theta}||^2 + \frac{1}{f_1}||\hat{h}^{i}_{\theta}||^2. 
\end{equation*}

\noindent Multiply \eqref{term_1} by $i^{2q}$ and sum over all $1 \leq i \leq K$. The result follows by collecting all the inequalities above. 

\bigskip

\begin{lemma}
We have the estimate
    \begin{equation*}
\partial_{t}\sum_{i=1}^{K}||i^{q}\nabla_{x}\hat{h}^{i}_{\theta}||^2 \leq -\frac{\lambda_{D}}{\varepsilon}\sum_{i=1}^{K} ||i^{q}\nabla_{x}\hat{h}^{i,\perp}_{\theta}||^{2}_{\Lambda} +R_2, 
\end{equation*}

\noindent where 

\begin{equation*}
    R_2 = f_2\sum_{i=1}^{K}||i^q \nabla_{x}A^{i}_{\theta}||^2 + \frac{1}{f_2}\sum_{i=1}^{K}||i^q\ \nabla_{x}\hat{h}^{i}_{\theta}||^2 + C\sum_{i=1}^{K}i^{2q}\mathcal{R}_{\theta,t}^{b,i,\nabla_{x}}+\sum_{i=1}^{K}i^{2q}r_{3,t}^{\nabla_x,i}, 
\end{equation*}

\noindent for $\lambda_{D}, f_2, C>0$, all independent of $\varepsilon$ with $f_2$ a free variable to be chosen later. 
\end{lemma}

\textit{Proof of lemma.} Take $\nabla_{x}$ on both sides of \eqref{Boltzmann_eq_SG_NN} and then take the inner product with respect to $\nabla_{x}\hat{h}^{i}_{\theta}$, we get
\begin{equation}\label{term_2}
    \begin{aligned}
        \ \frac{1}{2}\partial_{t}||\nabla_{x}\hat{h}^{i}_{\theta}||^2 = & -\langle v \cdot \nabla_{x}\nabla_{x}\hat{h}^{i}_{\theta}, \nabla_{x}\hat{h}^{i}_{\theta}\rangle + \frac{1}{\varepsilon}\sum_{k=1}^{K}\langle \mathcal{L}_{ik}(\nabla_{x}\hat{h}^{k}_{\theta}), \nabla_{x}\hat{h}^{i}_{\theta} \rangle + \langle \nabla_{x}A_{\theta}^{i}, \nabla_{x}\hat{h}^{i}_{\theta} \rangle. 
    \end{aligned}
\end{equation}

\noindent Following exactly the same arguments as for $||\hat{h}^{i}_{\theta}||^2$, we can derive 
\begin{equation*}
    \begin{aligned}
        \ &-\langle v \cdot \nabla_{x}\nabla_{x}\hat{h}^{i}_{\theta}, \nabla_{x}\hat{h}^{i}_{\theta}\rangle \leq C\mathcal{R}_{\theta,t}^{b,i,\nabla_{x}}+r_{3,t}^{\nabla_{x},i}\\[2pt]
        \ &\sum_{i=1}^{K}i^{2q}\sum_{k=1}^{K}\langle \mathcal{L}_{ik}(\nabla_{x}\hat{h}^{k}_{\theta}), \nabla_{x}\hat{h}^{i}_{\theta} \rangle \leq -\lambda_{D}\sum_{i=1}^{K} i^{2q}||\nabla_{x}\hat{h}^{i,\perp}_{\theta}||^{2}_{\Lambda},  
    \end{aligned}
\end{equation*}
\noindent and 
\begin{equation*}
    \langle \nabla_{x}A^{i}_{\theta}, \nabla_{x}\hat{h}^{i}_{\theta} \rangle \leq f_2||\nabla_{x}A^{i}_{\theta}||^2 + \frac{1}{f_2}||\nabla_{x}\hat{h}^{i}_{\theta}||^{2}. 
\end{equation*}

\noindent Hence the result follows by multiplying \eqref{term_2} by $i^2q$ and sum over all $1 \leq i \leq K$. 

\bigskip

\begin{lemma}
We have the estimate
    \begin{equation*}
    \begin{aligned}
        \ \partial_t \sum_{i=1}^{K}||i^{q}\nabla_{v}\hat{h}^{i,\perp}_{\theta}||^2 \leq & \frac{C_1}{\varepsilon}\sum_{i=1}^{K}||i^{q}\hat{h}^{k,\perp}_{\theta}||^2_{\Lambda}- \frac{C_2}{\varepsilon} \sum_{i=1}^{K}||i^{q}\nabla_{v}\hat{h}^{i,\perp}_{\theta}||^2_{\Lambda} \\
        \ & + \varepsilon C_3\sum_{i=1}^{K}||i^{q}\nabla_{x}\hat{h}^{i}_{\theta}||^2 +R_3, 
    \end{aligned}
\end{equation*}

\noindent where 
\begin{equation*}
    \begin{aligned}
        \ R_3 = & f_3\tilde{C}_{1}(\sum_{i=1}^{K}||i^{q}A^{i}_{\theta}||^2 + \sum_{i=1}^{K}||i^{q}\nabla_{v}A^{i}_{\theta}||^2) + \frac{1}{f_3}\tilde{C}_{1}(\sum_{i=1}^{K}||i^{q}\hat{h}^{i}_{\theta}||^2 + \sum_{i=1}^{K}||i^{q}\nabla_{v}\hat{h}^{i}_{\theta}||^2) \\[2pt]
        \ & + C_1^{'}(\sum_{i=1}^{K}i^{2q}\mathcal{R}_{\theta,t}^{b,i} + \sum_{i=1}^{K}i^{2q}\mathcal{R}_{\theta,t}^{b,i,\nabla_{v}}) + C_{2}^{'}(\sum_{i=1}^{K}i^{2q}r_{3,t}^{i} + \sum_{i=1}^{K}i^{2q}r_{3,t}^{\nabla_{v},i}), 
    \end{aligned}
\end{equation*}
\noindent for $C_1,C_2,C_3, C_1^{'},C_2^{'}, \tilde{C_1}, f_3 >0$, all independent of $\varepsilon$ with $f_3$ a free variable to be chosen later. 
\end{lemma}

\textit{Proof of lemma.} Take $I - \pi_{\mathcal{L}}$ and $\nabla_{v}$ respectively on both sides of \eqref{Boltzmann_eq_SG_NN}, and then take the inner product with respect to $ \nabla_{v}\hat{h}^{i,\perp}_{\theta}$. We get 
\begin{equation}\label{term_3}
    \begin{aligned}
        \ \frac{1}{2}\partial_{t}||\nabla_{v}\hat{h}^{i,\perp}_{\theta}||^2 =& -\langle \nabla_{v}(v \cdot \nabla_{x}\hat{h}^{i}_{\theta})^{\perp}, \nabla_{v}\hat{h}^{i,\perp}_{\theta}\rangle + \frac{1}{\varepsilon}\sum_{k=1}^{K}\langle \nabla_{v}\mathcal{L}_{ik}(\hat{h}^{k,\perp}_{\theta}), \nabla_{v}\hat{h}^{i,\perp}_{\theta} \rangle \\
        \ & + \langle \nabla_{v}A^{i,\perp}_{\theta}, \nabla_{v}\hat{h}^{i,\perp}_{\theta} \rangle, 
    \end{aligned}
\end{equation}
\noindent where we used $\mathcal{L}_{ik}(f) \in \mathcal{N}(\mathcal{L}_{ik})^{\perp} ~\forall f \in L_{v}^{2}$, and hence $(\mathcal{L}_{ik}(\hat{h}^{k}_{\theta}))^{\perp} = \mathcal{L}_{ik}(\hat{h}^{k}_{\theta}) - \pi_{\mathcal{L}}(\mathcal{L}_{ik}(\hat{h}^{k}_{\theta}))=\mathcal{L}_{ik}(\hat{h}^{k}_{\theta})=\mathcal{L}_{ik}(\hat{h}^{k,\perp}_{\theta})$. 

First, consider the term $\frac{1}{\varepsilon}\sum_{k=1}^{K}\langle \nabla_{v}\mathcal{L}_{ik}(\hat{h}^{k,\perp}_{\theta}), \nabla_{v}\hat{h}^{i,\perp}_{\theta} \rangle$. we multiply the term by $i^{2q}$ and sum over all $0 \leq i \leq K$. Then we use $\bm{H1}-\bm{H2}$ and follow exactly the same arguments in \cite{LiuJin2018} to get 
\begin{equation*}
    \begin{aligned}
        \ & \frac{1}{\varepsilon}\sum_{i=1}^{K}i^{2q}\sum_{k=1}^{K}\langle \nabla_{v}\mathcal{L}_{ik}(\hat{h}^{k,\perp}_{\theta}), \nabla_{v}\hat{h}^{i,\perp}_{\theta} \rangle = \frac{1}{\varepsilon}\sum_{i,k=1}^{K}i^{2q}\langle \nabla_{v}\mathcal{L}_{ik}(\hat{h}^{k,\perp}_{\theta}), \nabla_{v}\hat{h}^{i,\perp}_{\theta} \rangle \\
        \ & \leq \frac{1}{\varepsilon}\sum_{i,k=1}^{K}\chi_{ik} \cdot i^{2q}\bigl( (C(\delta)\frac{\nu_{1}^{\Lambda}}{\nu_{0}^{\Lambda}} + \nu^{\Lambda}_{4})||\hat{h}^{k,\perp}_{\theta}||^2_{\Lambda} + (\delta \frac{\nu_{1}^{\Lambda}}{\nu_{0}^{\Lambda}}-\nu^{\Lambda}_{3})||\nabla_{v}\hat{h}^{i,\perp}_{\theta}||^2_{\Lambda}\bigr) \\
        \ & = \frac{1}{\varepsilon}\sum_{i,k=1}^{K}\chi_{ik} \cdot \frac{i^{2q}}{k^{2q}}(C(\delta)\frac{\nu_{1}^{\Lambda}}{\nu_{0}^{\Lambda}} + \nu^{\Lambda}_{4})||k^{q}\hat{h}^{k,\perp}_{\theta}||^2_{\Lambda} + \frac{1}{\varepsilon}\sum_{i,k=1}^{K}\chi_{ik} \cdot (\delta \frac{\nu_{1}^{\Lambda}}{\nu_{0}^{\Lambda}}-\nu^{\Lambda}_{3})||i^{q}\nabla_{v}\hat{h}^{i,\perp}_{\theta}||^2_{\Lambda} \\
        \ & \leq \frac{3 \times 4^q}{\varepsilon}(C(\delta)\frac{\nu_{1}^{\Lambda}}{\nu_{0}^{\Lambda}} + \nu^{\Lambda}_{4})\sum_{k=1}^{K}||k^{q}\hat{h}^{k,\perp}_{\theta}||^2_{\Lambda} + \frac{3}{\varepsilon}(\delta \frac{\nu_{1}^{\Lambda}}{\nu_{0}^{\Lambda}}-\nu^{\Lambda}_{3})\sum_{i=1}^{K}||i^{q}\nabla_{v}\hat{h}^{i,\perp}_{\theta}||^2_{\Lambda}. 
    \end{aligned}
\end{equation*}

For the term $\ -\langle \nabla_{v}(v \cdot \nabla_{x}\hat{h}^{i}_{\theta})^{\perp}, \nabla_{v}\hat{h}^{i,\perp}_{\theta}\rangle$, we expand it as 
\begin{equation*}
    \begin{aligned}
        \ -\langle \nabla_{v}(v \cdot \nabla_{x}\hat{h}^{i}_{\theta})^{\perp}, \nabla_{v}\hat{h}^{i,\perp}_{\theta}\rangle =& -\langle \nabla_{v}(v \cdot \nabla_{x}\hat{h}^{i}_{\theta}), \nabla_{v}\hat{h}^{i,\perp}_{\theta}\rangle + \langle \nabla_{v}\pi_{\mathcal{L}}(v \cdot \nabla_{x}\hat{h}^{i}_{\theta}), \nabla_{v}\hat{h}^{i,\perp}_{\theta}\rangle \\[2pt]
        \ =& -\langle \nabla_{x}\hat{h}^{i}_{\theta},\nabla_{v}\hat{h}^{i,\perp}_{\theta} \rangle - \langle v \cdot \nabla_{v}\nabla_{x}\hat{h}^{i}_{\theta},\nabla_{v}\hat{h}^{i,\perp}_{\theta} \rangle \\[2pt]
        \ & + \langle \nabla_{v}\pi_{\mathcal{L}}(v \cdot \nabla_{x}\hat{h}^{i}_{\theta}),\nabla_{v}\hat{h}^{i,\perp}_{\theta} \rangle \\[2pt]
        \ =& -\langle \nabla_{x}\hat{h}^{i}_{\theta},\nabla_{v}\hat{h}^{i,\perp}_{\theta} \rangle - \langle v \cdot \nabla_{v}\nabla_{x}\pi_{\mathcal{L}}(\hat{h}^{i}_{\theta}),\nabla_{v}\hat{h}^{i,\perp}_{\theta} \rangle \\[2pt]
        \ & - \langle v \cdot \nabla_{x}\nabla_{v}\hat{h}^{i,\perp}_{\theta},\nabla_{v}\hat{h}^{i,\perp}_{\theta} \rangle + \langle \nabla_{v}\pi_{\mathcal{L}}(v \cdot \nabla_{x}\hat{h}^{i}_{\theta}),\nabla_{v}\hat{h}^{i,\perp}_{\theta} \rangle. 
    \end{aligned}
\end{equation*}

Using integration by part and inequality \eqref{H_3_1}, we have 
\begin{equation*}
    - \langle v \cdot \nabla_{x}\nabla_{v}\hat{h}^{i,\perp}_{\theta},\nabla_{v}\hat{h}^{i,\perp}_{\theta} \rangle \leq C_1^{'}(\mathcal{R}_{\theta,t}^{b,i} + \mathcal{R}_{\theta,t}^{b,i,\nabla_{v}}) + C_{2}^{'}(r_{3,t}^{i} + r_{3,t}^{\nabla_{v},i}). 
\end{equation*}

Apply the Cauchy-Schwarz, Young's inequality and inequality \eqref{H_1_1}, we have
\begin{equation*}
    \begin{aligned}
        \ -\langle \nabla_{x}\hat{h}^{i}_{\theta},\nabla_{v}\hat{h}^{i,\perp}_{\theta} \rangle \leq& \eta_{1}||\nabla_{x}\hat{h}^{i}_{\theta}||^2 + \frac{1}{\eta_{1}}||\nabla_{v}\hat{h}^{i,\perp}_{\theta}||^2 \\
        \ \leq & \eta_{1}||\nabla_{x}\hat{h}^{i}_{\theta}||^2 + \frac{\nu_{1}^{\Lambda}}{\eta_{1}\nu_{0}^{\Lambda}}||\nabla_{v}\hat{h}^{i,\perp}_{\theta}||^{2}_{\Lambda}. 
    \end{aligned}  
\end{equation*}

Apply the Cauchy-Schwarz, Young's inequality, inequality \eqref{H_3_1} and \eqref{H_1_1}, we can derive 
\begin{equation*}
    \begin{aligned}
        \ - \langle v \cdot \nabla_{v}\nabla_{x}\pi_{\mathcal{L}}(\hat{h}^{i}_{\theta}),\nabla_{v}\hat{h}^{i,\perp}_{\theta} \rangle \leq& \eta_{2}||v \cdot \nabla_{v}\nabla_{x}\pi_{\mathcal{L}}(\hat{h}^{i}_{\theta})||^2 + \frac{1}{\eta_{2}}||\nabla_{v}\hat{h}^{i,\perp}_{\theta}||^2\\[2pt]
        \ =& \eta_{2}||v \cdot \nabla_{v}\pi_{\mathcal{L}}(\nabla_{x}\hat{h}^{i}_{\theta})||^2 + \frac{1}{\eta_{2}}||\nabla_{v}\hat{h}^{i,\perp}_{\theta}||^2\\[2pt]
        \ \leq& \eta_{2}C_{\pi 1}||\pi_{\mathcal{L}}(\nabla_{x}\hat{h}^{i}_{\theta})||^2 + \frac{1}{\eta_{2}}||\nabla_{v}\hat{h}^{i,\perp}_{\theta}||^2 \\[2pt]
        \ \leq& \eta_{2}C_{\pi 1}||\nabla_{x}\hat{h}^{i}_{\theta}||^2 + \frac{\nu_{1}^{\Lambda}}{\eta_{2}\nu_{0}^{\Lambda}}||\nabla_{v}\hat{h}^{i,\perp}_{\theta}||^{2}_{\Lambda}
    \end{aligned}
\end{equation*}

\noindent and 
\begin{equation*}
    \begin{aligned}
        \ \langle \nabla_{v}\pi_{\mathcal{L}}(v \cdot \nabla_{x}\hat{h}^{i}_{\theta}),\nabla_{v}\hat{h}^{i,\perp}_{\theta} \rangle \leq& \eta_{3}||\nabla_{v}\pi_{\mathcal{L}}(v \cdot \nabla_{x}\hat{h}^{i}_{\theta})||^2 + \frac{1}{\eta_{3}}||\nabla_{v}\hat{h}^{i,\perp}_{\theta}||^2  \\[2pt]
        \ \leq& \eta_{3}C_{\pi 1}||\pi_{\mathcal{L}}(\nabla_{x}\hat{h}^{i}_{\theta})||^2 + \frac{1}{\eta_{3}}||\nabla_{v}\hat{h}^{i,\perp}_{\theta}||^2 \\[2pt]
        \ \leq& \eta_{3}C_{\pi 1}||\nabla_{x}\hat{h}^{i}_{\theta}||^2 + \frac{\nu_{1}^{\Lambda}}{\eta_{3}\nu_{0}^{\Lambda}}||\nabla_{v}\hat{h}^{i,\perp}_{\theta}||^{2}_{\Lambda}. 
    \end{aligned}
\end{equation*}

\noindent Finally, we use the Cauchy-Schwarz, Young's inequality and \eqref{H_3_1} to bound $\langle \nabla_{v}A^{i,\perp}_{\theta}, \nabla_{v}\hat{h}^{i,\perp}_{\theta} \rangle$ by
\begin{equation*}
    \begin{aligned}
        \ \langle \nabla_{v}A^{i,\perp}_{\theta}, \nabla_{v}\hat{h}^{i,\perp}_{\theta} \rangle \leq & f_{3}||\nabla_{v}A^{i,\perp}_{\theta}||^2 + \frac{1}{f_{3}}||\nabla_{v}\hat{h}^{i,\perp}_{\theta}||^{2} \\[2pt]
        \ \leq & f_3\tilde{C}_{1}(||A^{i}_{\theta}||^2 + ||\nabla_{v}A^{i}_{\theta}||^2) + \frac{1}{f_3}\tilde{C}_{1}(||\hat{h}^{i}_{\theta}||^2 + ||\nabla_{v}\hat{h}^{i}_{\theta}||^2). 
    \end{aligned}
\end{equation*}

Multiply both sides of \eqref{term_3} by $i^{2q}$ and sum over all $0 \leq i \leq K$. Then collect all the inequalities above, we get 
\begin{equation*}
    \begin{aligned}
        \ \partial_t \sum_{i=1}^{K}||i^{q}\nabla_{v}\hat{h}^{i,\perp}_{\theta}||^2 \leq & \frac{3 \times 4^q}{\varepsilon}(C(\delta)\frac{\nu_{1}^{\Lambda}}{\nu_{0}^{\Lambda}} + \nu^{\Lambda}_{4})\sum_{i=1}^{K}||i^{q}\hat{h}^{k,\perp}_{\theta}||^2_{\Lambda} \\[2pt]
        \ & + \bigl( \frac{3}{\varepsilon}(\delta \frac{\nu_{1}^{\Lambda}}{\nu_{0}^{\Lambda}}-\nu^{\Lambda}_{3}) + \frac{\nu_{1}^{\Lambda}}{\eta_{1}\nu_{0}} + \frac{\nu_{1}^{\Lambda}}{\eta_{2}\nu_{0}} + \frac{\nu_{1}^{\Lambda}}{\eta_{3}\nu_{0}} \bigr) \sum_{i=1}^{K}||i^{q}\nabla_{v}\hat{h}^{i,\perp}_{\theta}||^2_{\Lambda} \\[2pt]
        \ & + (\eta_1 + \eta_{2}C_{\pi 1} + \eta_{3}C_{\pi 1})\sum_{i=1}^{K}||i^{q}\nabla_{x}\hat{h}^{i}_{\theta}||^2 +R_3. 
    \end{aligned}
\end{equation*}

\noindent Let $\eta_1=\eta_2=\eta_3=\varepsilon \hat{\eta}$. We then choose $\hat{\eta}$ to be big enough and $\delta$ to be small enough to get the desired result. 

\bigskip

\begin{lemma}
We have the estimate
    \begin{equation*}
    \begin{aligned}
        \ \partial_{t}\sum_{i=1}^{K}i^{2q}\varepsilon\langle \nabla_{x}\hat{h}^{i}_{\theta}, \nabla_{v}\hat{h}^{i}_{\theta} \rangle  
        \ \leq & -\frac{\varepsilon}{4}\sum_{i=1}^{K}||i^{q}\nabla_{x}\hat{h}^{i}_{\theta}||^2 + \frac{D_1 e}{\varepsilon}\sum_{i=1}^{K} ||i^{q}\nabla_{x}\hat{h}^{i,\perp}_{\theta}||_{\Lambda}^{2} \\[2pt]
        \  & + \frac{D_2}{e\varepsilon}\sum_{i=1}^{K}||i^{q}\nabla_{v}\hat{h}^{i,\perp}_{\theta}||^{2}_{\Lambda} + R_4, 
    \end{aligned}
\end{equation*}

\noindent where  

\begin{equation*}
    \begin{aligned}
        R_4 = & f_4\sum_{i=1}^{K}||i^q\nabla_{x}A^{i}_{\theta}||^2 + \frac{1}{f_4}\sum_{i=1}^{K}||i^q\nabla_{v}\hat{h}^{i}_{\theta}||^2 \\[2pt]
        \ & + C_1^{'}(\sum_{i=1}^{K}i^{2q}\mathcal{R}_{\theta,t}^{b,i,\nabla_{x}} +  \sum_{i=1}^{K}i^{2q}\mathcal{R}_{\theta,t}^{b,i,\nabla_{v}}) + C_{2}^{'}(\sum_{i=1}^{K}i^{2q}r_{3,t}^{\nabla_{x},i} + \sum_{i=1}^{K}i^{2q}r_{3,t}^{\nabla_{v},i}), 
    \end{aligned}
\end{equation*}

\noindent for $D_1,D_2,f_4,e,C_1^{'}, C_2^{'}>0$, all independent of $\varepsilon$ with $f_4>0$ and $e>4$ two free variables to be chosen later. 
\end{lemma}

\textit{Proof of lemma.} Take $\nabla_{x}$ on both sides of \eqref{Boltzmann_eq_SG_NN} and then take the inner product with respect to $\nabla_{v}\hat{h}^{i}_{\theta}$, we get
\begin{equation}\label{term_4}
    \ \frac{1}{2}\partial_{t}\langle \nabla_{x}\hat{h}^{i}_{\theta}, \nabla_{v}\hat{h}^{i}_{\theta} \rangle = -\langle \nabla_{x}(v \cdot \nabla_{x}\hat{h}^{i}_{\theta}), \nabla_{v}\hat{h}^{i}_{\theta}\rangle + \frac{1}{\varepsilon}\sum_{k=1}^{K}\langle \mathcal{L}_{ik}(\nabla_{x}\hat{h}^{k}_{\theta}), \nabla_{v}\hat{h}^{i}_{\theta} \rangle + \langle \nabla_{x}A^{i}, \nabla_{v}\hat{h}^{i}_{\theta} \rangle. 
\end{equation}

Using integration by part, Cauchy-Schwarz and Young's inequality, one can show that
\begin{equation*}
    \begin{aligned}
        \ -\langle \nabla_{x}(v \cdot \nabla_{x}\hat{h}^{i}_{\theta}), \nabla_{v}\hat{h}^{i}_{\theta}\rangle \leq& \langle v \cdot \nabla_{x}\hat{h}^{i}_{\theta}, \nabla_{x}\nabla_{v}\hat{h}^{i}_{\theta}\rangle +C_1\mathcal{R}_{\theta,t}^{b,i,\nabla_{x}}+C_2\mathcal{R}_{\theta,t}^{b,i,\nabla_{v}} \\[2pt]
        \ & +C_{1}^{'}r_{3,t}^{\nabla_{x},i} +C_{2}^{'}r_{3,t}^{\nabla_{v},i} \\[2pt]
        \ =& -\langle \nabla_{v}(v \cdot \nabla_{x}\hat{h}^{i}_{\theta}), \nabla_{x}\hat{h}^{i}_{\theta}\rangle +C_1\mathcal{R}_{\theta,t}^{b,i,\nabla_{x}}+C_2\mathcal{R}_{\theta,t}^{b,i,\nabla_{v}} \\[2pt]
        \ & +C_{1}^{'}r_{3,t}^{\nabla_{x},i} +C_{2}^{'}r_{3,t}^{\nabla_{v},i} \\[2pt]
    \ =&-||\nabla_{x}\hat{h}^{i}_{\theta}||^2 - \langle v \cdot \nabla_{v}\nabla_{x}\hat{h}^{i}_{\theta}, \nabla_{x}\hat{h}^{i}_{\theta}\rangle +C_1\mathcal{R}_{\theta,t}^{b,i,\nabla_{x}} \\[2pt]
    \ & +C_2\mathcal{R}_{\theta,t}^{b,i,\nabla_{v}}+C_{1}^{'}r_{3,t}^{\nabla_{x},i} +C_{2}^{'}r_{3,t}^{\nabla_{v},i} \\[2pt]
    \ \leq & -||\nabla_{x}\hat{h}^{i}_{\theta}||^2 + \langle v \cdot \nabla_{v}\hat{h}^{i}_{\theta}, \nabla_{x}\nabla_{x}\hat{h}^{i}_{\theta}\rangle +2C_1\mathcal{R}_{\theta,t}^{b,i,\nabla_{x}} \\[2pt]
    \ & +2C_2\mathcal{R}_{\theta,t}^{b,i,\nabla_{v}}+2C_{1}^{'}r_{3,t}^{\nabla_{x},i} +2C_{2}^{'}r_{3,t}^{\nabla_{v},i}. 
    \end{aligned}
\end{equation*}

Since $\langle v \cdot \nabla_{v}\hat{h}^{i}_{\theta}, \nabla_{x}\nabla_{x}\hat{h}^{i}_{\theta}\rangle = \langle \nabla_{x}(v \cdot \nabla_{x}\hat{h}^{i}_{\theta}), \nabla_{v}\hat{h}^{i}_{\theta}\rangle$, we have 
\begin{equation*}
    \ -\langle \nabla_{x}(v \cdot \nabla_{x}\hat{h}^{i}_{\theta}), \nabla_{v}\hat{h}^{i}_{\theta}\rangle \leq -\frac{1}{2}||\nabla_{x}\hat{h}^{i}_{\theta}||^2  +C_1\mathcal{R}_{\theta,t}^{b,i,\nabla_{x}}+C_2\mathcal{R}_{\theta,t}^{b,i,\nabla_{v}}+C_{1}^{'}r_{3,t}^{\nabla_{x},i} +C_{2}^{'}r_{3,t}^{\nabla_{v},i}. 
\end{equation*}

By \eqref{H_1_3}, for each $\mathcal{L}_{ik}$, there is a constant $C^{\mathcal{L}_{ik}}$ such that $\langle \mathcal{L}_{ik}(h),g \rangle \leq C^{\mathcal{L}_{ik}}||h||_{\Lambda_{v}}||g||_{\Lambda_{v}} ~\forall h,g \in L_{v}^{2}$. Let $C^{\mathcal{L}}=\max\{ C^{\mathcal{L}_{ik}} \}_{0 \leq i,k \leq K}$. Then 
\begin{equation*}
    \begin{aligned}
        \ \frac{1}{\varepsilon}\sum_{k=1}^{K}\langle \mathcal{L}_{ik}(\nabla_{x}\hat{h}^{k}_{\theta}), \nabla_{v}\hat{h}^{i}_{\theta} \rangle =& \frac{1}{\varepsilon}\sum_{k=1}^{K}\langle \mathcal{L}_{ik}(\nabla_{x}\hat{h}^{k,\perp}_{\theta}), \nabla_{v}\hat{h}^{i}_{\theta} \rangle \\[2pt]
        \ \leq & \frac{1}{\varepsilon}\sum_{k=1}^{K} C^{\mathcal{L}}||\nabla_{x}\hat{h}^{k,\perp}_{\theta}||_{\Lambda}||\nabla_{v}\hat{h}^{i}_{\theta}||_{\Lambda}\cdot \chi_{ik} \\[2pt]
        \ \leq & \frac{C^{\mathcal{L}}\eta}{\varepsilon}\sum_{k=1}^{K} ||\nabla_{x}\hat{h}^{k,\perp}_{\theta}||_{\Lambda}^{2}\cdot \chi_{ik} + \frac{C^{\mathcal{L}}(K+1)}{\varepsilon \eta}||\nabla_{v}\hat{h}^{i}_{\theta}||^{2}_{\Lambda}
    \end{aligned}
\end{equation*}
 where we apply Young's inequality in the last inequality above. Using inequality \eqref{H_3_2} and the Poincare inequality \eqref{Poincare inequality}, We can further bound $||\nabla_{v}\hat{h}^{i}_{\theta}||^{2}_{\Lambda}$ by 
\begin{equation*}
    \begin{aligned}
        \ ||\nabla_{v}\hat{h}^{i}_{\theta}||^{2}_{\Lambda} =& ||\nabla_{v} \pi_{\mathcal{L}}(\hat{h}^{i}_{\theta}) + \nabla_{v}\hat{h}^{i,\perp}_{\theta}||^{2}_{\Lambda}\\[2pt]
        \ \leq& 2||\nabla_{v} \pi_{\mathcal{L}}(\hat{h}^{i}_{\theta})||^{2}_{\Lambda} +2||\nabla_{v}\hat{h}^{i,\perp}_{\theta}||^{2}_{\Lambda} \\[2pt]
        \ \leq& 2C_{\pi}||\pi_{\mathcal{L}}(\hat{h}^{i}_{\theta})||^2 + 2||\nabla_{v}\hat{h}^{i,\perp}_{\theta}||^{2}_{\Lambda}\\[2pt]
        \ \leq& 2C_{\pi}C_{p}||\nabla_{x}\hat{h}^{i}_{\theta}||^2 + 2||\nabla_{v}\hat{h}^{i,\perp}_{\theta}||^{2}_{\Lambda} 
    \end{aligned}
\end{equation*}

\noindent Finally, we use Cauchy-Schwarz, Young's inequality to bound $\langle \nabla_{x}A^{i}_{\theta}, \nabla_{v}\hat{h}^{i}_{\theta} \rangle$ by 

\begin{equation*}
         \langle \nabla_{x}A^{i}_{\theta}, \nabla_{v}\hat{h}^{i}_{\theta} \rangle \leq f_{4}||\nabla_{x}A^{i}_{\theta}||^2 + \frac{1}{f_{4}}||\nabla_{v}\hat{h}^{i}_{\theta}||^2. 
\end{equation*}

Set $\eta=\frac{D}{\varepsilon}$ for $D>0$ and collect all the inequalities above, we get
\begin{equation*}
    \begin{aligned}
        \ \partial_{t}\langle \nabla_{x}\hat{h}^{i}_{\theta}, \nabla_{v}\hat{h}^{i}_{\theta} \rangle \leq & \bigl( 2C_{\pi}C_{p}\frac{C^{\mathcal{L}}(K+1)}{D}-\frac{1}{2} \bigr)||\nabla_{x}\hat{h}^{i}_{\theta}||^2 + \frac{C^{\mathcal{L}}D}{\varepsilon^2}\sum_{k=1}^{K} ||\nabla_{x}\hat{h}^{k,\perp}_{\theta}||_{\Lambda}^{2}\cdot \chi_{ik} \\[2pt]
        \  & + \frac{2C^{\mathcal{L}}(K+1)}{D}||\nabla_{v}\hat{h}^{i,\perp}_{\theta}||^{2}_{\Lambda} + f_4||\nabla_{x}A^{i}_{\theta}||^2 + \frac{1}{f_4}||\nabla_{v}\hat{h}^{i}_{\theta}||^2 \\[2pt]
        \ & + C_1\mathcal{R}_{\theta,t}^{b,i,\nabla_{x}}+C_2\mathcal{R}_{\theta,t}^{b,i,\nabla_{v}}+C_{1}^{'}r_{3,t}^{\nabla_{x},i} +C_{2}^{'}r_{3,t}^{\nabla_{v},i}. 
    \end{aligned}
\end{equation*}

We set $D=2C_{\pi}C_{p}C^{\mathcal{L}}(K+1)e$ for some $e>4$ to be chosen later. Then we multiply both sides of the inequality above by $\varepsilon$ and use the fact that $\varepsilon <1$ to get  
\begin{equation}\label{bound for gradient_x_v intermediate}
    \begin{aligned}
        \ \partial_{t}\varepsilon\langle \nabla_{x}\hat{h}^{i}_{\theta}, \nabla_{v}\hat{h}^{i}_{\theta} \rangle \leq & -\frac{\varepsilon}{4}||\nabla_{x}\hat{h}^{i}_{\theta}||^2 + \frac{\Tilde{D_1} e}{\varepsilon}\sum_{k=1}^{K} ||\nabla_{x}\hat{h}^{k,\perp}_{\theta}||_{\Lambda}^{2}\cdot \chi_{ik} \\[2pt]
        \  & + \frac{D_2}{e\varepsilon}||\nabla_{v}\hat{h}^{i,\perp}_{\theta}||^{2}_{\Lambda} + f_4||\nabla_{x}A^{i}_{\theta}||^2 + \frac{1}{f_4}||\nabla_{v}\hat{h}^{i}_{\theta}||^2 \\[2pt]
        \ & + C_1\mathcal{R}_{\theta,t}^{b,i,\nabla_{x}}+C_2\mathcal{R}_{\theta,t}^{b,i,\nabla_{v}}+C_{1}^{'}r_{3,t}^{\nabla_{x},i} +C_{2}^{'}r_{3,t}^{\nabla_{v},i}. 
    \end{aligned}
\end{equation}

Finally, multiple both sides of \eqref{bound for gradient_x_v intermediate} by $i^{2q}$ and sum over all $1 \leq i \leq K$, we get 

\begin{equation*}
\begin{aligned}
        \ \partial_{t}\sum_{i=1}^{K}i^{2q}\varepsilon\langle \nabla_{x}\hat{h}^{i}_{\theta}, \nabla_{v}\hat{h}^{i}_{\theta} \rangle \leq & -\frac{\varepsilon}{4}\sum_{i=1}^{K}||i^{q}\nabla_{x}\hat{h}^{i}_{\theta}||^2 + \frac{\Tilde{D_1} e}{\varepsilon}\sum_{i,k=1}^{K} \frac{i^{2q}}{k^{2q}}||k^{q}\nabla_{x}\hat{h}^{k,\perp}_{\theta}||_{\Lambda}^{2}\cdot \chi_{ik} \\[2pt]
        \  & + \frac{D_2}{e\varepsilon}\sum_{i=1}^{K}||i^{q}\nabla_{v}\hat{h}^{i,\perp}_{\theta}||^{2}_{\Lambda} + R_4 \\[2pt]
        \ \leq & -\frac{\varepsilon}{4}\sum_{i=1}^{K}||i^{q}\nabla_{x}\hat{h}^{i}_{\theta}||^2 + \frac{D_1 e}{\varepsilon}\sum_{i=1}^{K} ||i^{q}\nabla_{x}\hat{h}^{i,\perp}_{\theta}||_{\Lambda}^{2} \\[2pt]
        \  & + \frac{D_2}{e\varepsilon}\sum_{i=1}^{K}||i^{q}\nabla_{v}\hat{h}^{i,\perp}_{\theta}||^{2}_{\Lambda} + R_4, 
    \end{aligned}
\end{equation*}

\noindent where $D_1 = 3 \times 4^q \Tilde{D_1}$ and the result follows. 

\bigskip
Collecting lemmas 5, 6, 7 and 8, we can derive that
\begin{equation}\label{bound for E_perp}
    \begin{aligned}
        \ \partial_{t}E_{\varepsilon \perp,t}^{K}(\hat{\bm{h}}_{\theta}) \leq &\frac{a_3C_1-a_1\lambda_D}{\varepsilon}\sum_{i=1}^{K} ||i^{q}\hat{h}^{i,\perp}_{\theta}||^{2}_{\Lambda} + \frac{a_4D_1e - a_2\lambda_D}{\varepsilon}\sum_{i=1}^{K} ||i^{q}\nabla_{x}\hat{h}^{i,\perp}_{\theta}||^{2}_{\Lambda} \\[2pt]
        \ &+\frac{\frac{a_4D_2}{e}-a_3C_2}{\varepsilon}\sum_{i=1}^{K}||i^{q}\nabla_{v}\hat{h}^{i,\perp}_{\theta}||^2_{\Lambda} + \varepsilon(a_3C_3 - \frac{a_4}{4})\sum_{i=1}^{K}||i^{q}\nabla_{x}\hat{h}^{i}_{\theta}||^2 + R, 
    \end{aligned}
\end{equation}
where $R=R_1 + R_2 + R_3 + R_4$. 

We first pick $a_3, a_4$ so that $a_3C_3 - \frac{a_4}{4}<0$. Then we pick $e$ big enough so that $\frac{a_4D_2}{e}-a_3C_2<0$. Then we pick $a_1$ big enough so that $a_3C_1-a_1\lambda_D<0$. Finally, we choose $a_2$ big enough so that $a_4D_1e - a_2\lambda_D<0$ and $E_{\varepsilon \perp,t}^{K}$ is equivalent to $E_{t}^{K}$. Then \eqref{bound for E_perp} becomes 
\begin{equation*}
    \begin{aligned}
        \ \partial_{t}E_{\varepsilon \perp,t}^{K}(\hat{\bm{h}}_{\theta}) \leq &-\frac{K_1}{\varepsilon}\sum_{i=1}^{K} ||i^{q}\hat{h}^{i,\perp}_{\theta}||^{2}_{\Lambda} - \frac{K_2}{\varepsilon}\sum_{i=1}^{K} ||i^{q}\nabla_{x}\hat{h}^{i,\perp}_{\theta}||^{2}_{\Lambda} \\[2pt]
        \ &- \frac{K_3}{\varepsilon}\sum_{i=1}^{K}||i^{q}\nabla_{v}\hat{h}^{i,\perp}_{\theta}||^2_{\Lambda} - \varepsilon K_4\sum_{i=1}^{K}||i^{q}\nabla_{x}\hat{h}^{i}_{\theta}||^2 +R\\[2pt]
        \ \leq &-\varepsilon \tilde{K} (\frac{\sum_{i=1}^{K} ||i^{q}\hat{h}^{i,\perp}_{\theta}||^{2}_{\Lambda}}{\varepsilon^2}+\frac{\sum_{i=1}^{K} ||i^{q}\nabla_{x}\hat{h}^{i,\perp}_{\theta}||^{2}_{\Lambda}}{\varepsilon^2} \\[2pt]
        \ &~~~~~~~~~+ \frac{\sum_{i=1}^{K}||i^{q}\nabla_{v}\hat{h}^{i,\perp}_{\theta}||^2_{\Lambda}}{\varepsilon^2} + \sum_{i=1}^{K}||i^{q}\nabla_{x}\hat{h}^{i}_{\theta}||^2) +R\\[2pt]
        \ &\leq -\varepsilon \tilde{K}(\sum_{i=1}^{K} ||i^{q}\hat{h}^{i,\perp}_{\theta}||^{2}_{\Lambda} + \sum_{i=1}^{K} ||i^{q}\nabla_{x}\hat{h}^{i,\perp}_{\theta}||^{2}_{\Lambda} + \sum_{i=1}^{K}||i^{q}\nabla_{v}\hat{h}^{i,\perp}_{\theta}||^2_{\Lambda} + \sum_{i=1}^{K}||i^{q}\nabla_{x}\hat{h}^{i}_{\theta}||^2) +R\\[2pt]
        \ & \leq -\varepsilon \tilde{K}E_{t}^{K}(\hat{\bm{h}}_{\theta}) + R. 
    \end{aligned}
\end{equation*}

\noindent The last inequality is based on the following fact: Let 

\begin{equation*}
    \begin{cases}
        \ F_1(h) = ||h^{\perp}||^2_{\Lambda}+ ||\nabla_{x}h^{\perp}||^2_{\Lambda} + ||\nabla_{v}h^{\perp}||^2_{\Lambda} + ||\nabla_{x}h||^2 \\[2pt]
        \ F_2(h) = ||h||^2_{\Lambda}+ ||\nabla_{x}h||^2_{\Lambda} + ||\nabla_{v}h^{\perp}||^2_{\Lambda}, 
    \end{cases}
\end{equation*}
then use \eqref{H_3_2} and Poincare inequality \eqref{Poincare inequality}, we have 
\begin{equation*}
    ||h||^2_{\Lambda} \leq ||h^{\perp}||^2_{\Lambda} + ||\pi_{\mathcal{L}}(h)||^2_{\Lambda} \leq ||h^{\perp}||^2_{\Lambda} + C_{\pi}C_{p}||\nabla_{x}h||^2. 
\end{equation*} 

\noindent Similarly, 
\begin{equation*}
    ||\nabla_{x}h||^2_{\Lambda} \leq ||\nabla_{x}h^{\perp}||^2_{\Lambda} + C_{\pi}||\nabla_{x}h||^2, 
\end{equation*}
\noindent so $F_2(h) \leq CF_{1}(h)$ for some $C>0$. On the other hand, using \eqref{H_3_2} and \eqref{H_1_1}, we get 
\begin{equation*}
    ||\nabla_{v}h||^2_{\Lambda} \leq ||\nabla_{v}h^{\perp}||^2_{\Lambda} + ||\nabla_{v}\pi_{\mathcal{L}}(h)||^2_{\Lambda} \leq ||\nabla_{v}h^{\perp}||^2_{\Lambda} + \frac{C_{\pi}\nu^{\Lambda}_{1}}{\nu^{\Lambda}_{0}}||h||^2_{\Lambda}. 
\end{equation*}
\noindent Hence $||h||_{H_{\Lambda}^{1}}^{2} \leq C^{'}F_{2}(h)$ for some $C^{'}>0$. Finally, since $\Lambda$-norm controls the standard $L^2$ norm, we have $ ||h||_{H_{x,v}^{1}}^{2} \leq C^{''}||h||_{H_{\Lambda}^{1}}^{2}$ for some $C^{''}>0$. Combining all these inequalities, we get $||h||_{H_{x,v}^{1}}^{2} \leq \tilde{C}F_{1}(h)$ for some $\tilde{C}>0$. 

\bigskip

Finally, we analyze the term $R$. Set $f_1 = f_2 = f_3 = f_4 = \frac{1}{f\varepsilon}$ for some $f>0$, we have  
\begin{equation*}
    R \leq \varepsilon fC_1E_{t}^{K}(\hat{\bm{h}}_{\theta}) + \frac{C_1}{\varepsilon f}E_{t}^{K}(\bm{A}_{\theta}) + C_2(\mathcal{R}_{\theta,t,sto}^{b} + \mathcal{R}_{\theta,t,sto}^{b,\nabla_{x}} + \mathcal{R}_{\theta,t,sto}^{b,\nabla_{v}}) + C_3(r_{3,t} + r_{3,t}^{\nabla_{x}} + r_{3,t}^{\nabla_{v}}), 
\end{equation*}

\noindent for some $C_1, C_2, C_3 >0$. We choose $f$ to be small enough so that $fC_1 - \tilde{K} <0$. Then 
\begin{equation*}
    \begin{aligned}
        \partial_{t}E_{\varepsilon \perp,t}^{K}(\hat{\bm{h}}_{\theta}) \leq & -\varepsilon K^{*}E_{t}^{K}(\hat{\bm{h}}_{\theta}) + \frac{C}{\varepsilon}(E_{t}^{K}(\bm{A}_{\theta})+\mathcal{R}_{\theta,t,sto}^{b} + \mathcal{R}_{\theta,t,sto}^{b,\nabla_{x}} + \mathcal{R}_{\theta,t,sto}^{b,\nabla_{v}} \\[2pt]
        \ & + r_{3,t} + r_{3,t}^{\nabla_{x}} + r_{3,t}^{\nabla_{v}}). 
    \end{aligned}
\end{equation*}

Split $\int_{\mathbb{R}^3} = \int_{\mathcal{D}} + \int_{\mathbb{R}^3-\mathcal{D}}$ and use assumption 2, we have
\begin{equation*}
    \begin{aligned}
        \ ||A^{i}_{\theta}||^2 \leq& ||\bm{d}^{1,i}_{\theta} \cdot \bm{\varphi}M||^2 +||d^{2,i}_{\theta}||^2 \\[2pt]
        \ =& \mathcal{R}_{\theta,t}^{1,i} + \mathcal{R}_{\theta,t}^{2,i} + r_{1,t}^{i}
    \end{aligned}
\end{equation*}

\noindent Similarly, $||\nabla_{x}A^{i}_{\theta}||^2 \leq \mathcal{R}_{\theta,t}^{1,i,\nabla_{x}} + \mathcal{R}_{\theta,t}^{2,i,\nabla_{x}} + r_{1,t}^{\nabla_{x},i}$ and $||\nabla_{v}A^{i}_{\theta}||^2 \leq \mathcal{R}_{\theta,t}^{1,i,\nabla_{v}} + \mathcal{R}_{\theta,t}^{2,i,\nabla_{v}} + r_{1,t}^{\nabla_{v},i}$. So we get 
\begin{equation*}
    \partial_{t}E_{\varepsilon \perp,t}^{K}(\hat{\bm{h}}_{\theta}) \leq -\varepsilon K^{*}E_{t}^{K}(\hat{\bm{h}}_{\theta}) + \frac{C}{\varepsilon}\tilde{R}
\end{equation*}
 where 
\begin{equation*}
    \begin{aligned}
        \ \tilde{R} = & \mathcal{R}_{\theta,t,sto}^{1} + \mathcal{R}_{\theta,t,sto}^{1,\nabla_{x}} + \mathcal{R}_{\theta,t,sto}^{1,\nabla_{v}}+ \mathcal{R}_{\theta,t,sto}^{2} + \mathcal{R}_{\theta,t,sto}^{2,\nabla_{x}} + \mathcal{R}_{\theta,t,sto}^{2,\nabla_{v}} +\mathcal{R}_{\theta,t,sto}^{b} \\[2pt]
        \ & + \mathcal{R}_{\theta,t,sto}^{b,\nabla_{x}} + \mathcal{R}_{\theta,t,sto}^{b,\nabla_{v}} + r_{1,t} + r_{1,t}^{\nabla_{x}} + r_{1,t}^{\nabla_{v}} + r_{3,t} + r_{3,t}^{\nabla_{x}} + r_{3,t}^{\nabla_{v}}
    \end{aligned}
\end{equation*}

\noindent Thus by Gronwall's inequality, we have 
\begin{equation*}
    \begin{aligned}
       \ E_{t}^{K}(\hat{\bm{h}}_{\theta}) \leq & \left(\frac{C}{\varepsilon}\int_{0}^{T} \tilde{R} ~dt + E_{0}^{K}(\hat{\bm{h}}_{\theta})\right)e^{-\varepsilon \tau t} \\[2pt]
       \ \leq & \frac{\tilde{C}\left(\mathcal{R}_{\theta,sto}^{G,H^1} + \delta(\varepsilon)\right)}{\varepsilon}e^{-\varepsilon \tau t}, 
    \end{aligned}
\end{equation*}
for some $\tilde{C}, \tau >0$ independent of $\varepsilon$. This completes the proof.

\begin{remark}
    Use similar arguments and a Gronwall type inequality, we can also show that $E_{t}^{K}(\hat{\bm{h}}_{\theta}) \leq \tilde{C}\left(\mathcal{R}_{\theta,sto}^{G,H^1} + \delta(\varepsilon)\right)$ for some $\tilde{C}>0$ independent of $\varepsilon$, hence the error won't explode when $\varepsilon <<1$. But in this case we can't derive the exponential decay in time. This is not a real problem since $e^{-\varepsilon \tau t} \rightarrow 1$ as $\varepsilon 
 \rightarrow 0$, hence the exponential decay becomes insignificant when we are in the fluid regime. We comment that in the acoustic scaling, we can't expect an exponential decay with decay rate independent of $\varepsilon$, as suggested in \cite{LiuJin2018}. 
\end{remark} 
\section{Conclusion}

In this paper, we construct an asymptotic-preserving neural network (APNN) for the linearized Boltzmann equation in the acoustic scaling with uncertainties. We employ the micro-macro decomposition method and build the APNN loss function based on the micro-macro system. Rigorous analysis have been conducted to show the existence of neural networks when the APNN loss goes to zero, and the convergence of the neural network approximated solution when the loss function tends to zero. In the future work, we will further improve our analysis result by incorporating the Barron-type functions \cite{GN2022} and develop posterior estimates for the neural network approximations.

\section*{Acknowledgement}
L. Liu acknowledges the support by National Key R\&D Program of China (2021YFA1001200), Ministry of Science and Technology in China, Early Career Scheme (24301021) and General Research Fund (14303022 \& 14301423) funded by Research Grants Council of Hong Kong.

\bibliographystyle{spmpsci}
\bibliography{Ref.bib}

\end{document}